\documentclass[a4paper, 10pt, parskip=half, english]{article}

\usepackage{babel}
\usepackage[T1]{fontenc}
\usepackage[utf8]{inputenc}
\usepackage{lmodern}
\usepackage{amsmath,bm,mathtools,amssymb}
\usepackage{graphicx,color}
\usepackage[footnotesize]{caption}
\usepackage{booktabs}
\usepackage{multirow}
\usepackage{relsize}

\usepackage[a4paper, margin=3cm]{geometry}
\usepackage[square,sort&compress,comma,numbers]{natbib}

\usepackage[justification=justified,singlelinecheck=false]{subcaption}

\usepackage{siunitx}
\usepackage{hyperref}

\date{\today}

\newcommand{\mR}{\mathbb{R}}

\newcommand{\GL}{G\!L}

\newcommand{\SO}{S\!O}
\let\so\undefined
\newcommand{\so}{\mathfrak{so}}
\newcommand{\SE}{{S\!E}}
\newcommand{\se}{\mathfrak{se}}

\newcommand{\Skw}{\operatorname{Skw}}

\newcommand{\Exp}{\operatorname{Exp}}
\newcommand{\Log}{\operatorname{Log}}

\newcommand{\tr}{\operatorname{tr}}

\newcommand{\T}{^{\!\mathrm{T}}}

\newcommand{\diff}[1][]{\mathrm{d}#1}

\newcommand{\pd}[2]{\frac{\partial #1}{\partial #2}}

\newcommand{\vga}{\bm{\gamma}}

\newcommand{\vth}{\bm{\theta}}

\newcommand{\vka}{\bm{\kappa}}

\newcommand{\vph}{\bm{\phi}} 

\newcommand{\vps}{\bm{\psi}}
\newcommand{\vom}{\bm{\omega}}

\newcommand{\vvep}{\bm{\varepsilon}}

\newcommand{\vOm}{\bm \Omega}

\newcommand{\va}{\mathbf a}
\newcommand{\vb}{\mathbf b}
\newcommand{\vc}{\mathbf c}
\newcommand{\vd}{\mathbf d}
\newcommand{\ve}{\mathbf e}
\newcommand{\vf}{\mathbf f}

\newcommand{\vm}{\mathbf m}
\newcommand{\vn}{\mathbf n}

\newcommand{\vq}{\mathbf q}
\newcommand{\vr}{\mathbf r}
\newcommand{\vs}{\mathbf s}

\newcommand{\vu}{\mathbf u}
\newcommand{\vv}{\mathbf v}

\newcommand{\vA}{\mathbf A}
\newcommand{\vB}{\mathbf B}
\newcommand{\vC}{\mathbf C}

\newcommand{\vH}{\mathbf H}
\newcommand{\vI}{\mathbf I}

\newcommand{\vM}{\mathbf M}

\newcommand{\vT}{\mathbf T}

\newcommand{\cE}{\mathcal{E}}

\newcommand{\cI}{\mathcal{I}}

\newcommand{\cK}{\mathcal{K}}

\title{
	\bfseries A family of total Lagrangian Petrov--Galerkin Cosserat rod finite element formulations
}

\author{
	Simon R. Eugster$^1$  and Jonas Harsch$^1$ \\[0.25cm]
	$^1$Institute for Nonlinear Mechanics \\
	University of Stuttgart \\
	Pfaffenwaldring 9, 70569 Stuttgart, Germany
}

\begin{document}
\maketitle
\begin{abstract}
The standard in rod finite element formulations is the Bubnov--Galerkin projection method, where the test functions arise from a consistent variation of the ansatz functions. This approach becomes increasingly complex when highly nonlinear ansatz functions are chosen to approximate the rod's centerline and cross-section orientations. Using a Petrov--Galerkin projection method, we propose a whole family of rod finite element formulations where the nodal generalized virtual displacements and generalized velocities are interpolated instead of using the consistent variations and time derivatives of the ansatz functions. This approach leads to a significant simplification of the expressions in the discrete virtual work functionals. In addition, independent strategies can be chosen for interpolating the nodal centerline points and cross-section orientations. We discuss three objective interpolation strategies and give an in-depth analysis concerning locking and convergence behavior for the whole family of rod finite element formulations.
\end{abstract}
\section{Introduction}
The theory of shear-deformable (spatial) rod formulations dates back to the pioneering works of the Cosserat brothers \cite{Cosserat1909}, Timoshenko \cite{Timoshenko1921}, Reissner \cite{Reissner1981} and Simo \cite{Simo1985}. Thus, depending on the chosen literature, shear-deformable rods are called among others (special) Cosserat rods \cite{Antman2005}, Simo--Reissner beams \cite{Meier2019}, spatial Timoshenko beams \cite{Eugster2015a}, or geometrically exact beams \cite{Betsch2002}. In this paper, we call the object of interest Cosserat rod, or just rod. 

For the numerical treatment of Cosserat rods, there is a vast amount of finite element formulations available, see \cite{Meier2019} for an exhaustive literature survey. Instead of giving an incomplete list of previous formulations, we want to address here some major challenges that appear in large-strain spatial rod finite element formulations. This could be of particular interest for application-oriented researchers, who are looking for the most efficient way to learn about rod theory and its computational analysis.

There have been many efforts to deduce the rod equations from a three-dimensional continuum theory, \cite{Antman1976a, Parker1979a, Auricchio2008, Eugster2015a}. However, the most direct way is to consider a rod in the sense of an intrinsic theory as a generalized one-dimensional continuum, i.e., as a spatial curve with additional degrees of freedom that describe the orientations of the rod's cross-sections. The procedure to obtain the governing rod equations in the form of partial differential equations (PDEs) can be sketched as follows, consult \cite{Eugster2020b, dellIsola2022b} for more details. First, one must accept to postulate mechanics within the variational framework of the principle of virtual work, see \cite{dellIsola2020b}. Second, the internal virtual work functional is defined. The structure of this functional follows from an objectivity requirement to the strain energy density. Then the internal virtual work functional together with the principle of virtual work implies the applicable external interactions a rod can resist to. For a Cosserat rod, these are forces and moments. Augmented by an inertial virtual work functional, finally, the partial differential equations of the Cosserat rod follow from the principle of virtual work, which can be seen as the weak form of the governing PDEs. The PDEs can also be interpreted as localized balances of linear and angular momentum. Hence, it is also possible to derive the governing PDEs from these two balance laws together with the assumption that the contact interactions inside the rod are given by forces and moments, \cite[Chapter 8]{Antman2005}. However, we advocate for the variational approach because the idea of finite elements is the approximation of the weak form of the governing PDEs, i.e., the approximation of the principle of virtual work. 

Besides the two different postulation philosophies in mechanics, another difficulty arises in the representation of vector quantities with respect to different bases. Since Truesdell and Noll \cite{Truesdell1965}, most theoretical works in continuum mechanics abstain from choosing a particular basis for the abstract three-dimensional real inner-product space, which in classical mechanics often serves as a convenient model of the ambient space. Only after the choice of a particular basis, a vector representation with three real numbers is obtained. We will discuss this issue in greater detail in Section \ref{sec:Model_of_the_ambient_space}. In this so-called coordinate-free framework, a distinction can be made whether the continuum is formulated in the material or spatial description, i.e., whether the appearing fields are formulated as functions of a reference configuration or of the deformed configuration, \cite{Eugster2022b, dellIsola2022b}. Partially, this concept can also be transferred to the theory of rods, see \cite{Simo1985}. However, in the reference configuration of a rod with a straight undeformed state, one usually introduces a distinct orthonormal basis, where two referential base vectors are aligned with the rod's cross-section-fixed bases \cite{Simo1985, Cardona1988, Gosh2008a}. In the subsequent computational analysis, it is then tacitly assumed to make use of this referential basis although another choice would have been possible. The components of material vectors with respect to the referential basis then coincide with the corresponding spatial vectors represented with respect to the cross-section-fixed bases. Consequently, it is enough to focus on the spatial description. Representations of spatial vectors with respect to cross-section-fixed bases or the inertial basis coincide with what in computational mechanics literature is referred to as the material and spatial description, respectively.

Maybe the major challenge in large-strain Cosserat rod finite element formulations is caused by the description of the cross-section orientations, which mathematically are captured by a function whose codomain is the set of orthogonal matrices, i.e., the special orthogonal group $\SO(3)$. This group is also a smooth three-dimensional manifold and infinitely many different $\SO(3)$-parametrizations exist, e.g., Euler-angles, Tait-Bryan angles, rotation vectors, Rodrigues parameters. Unfortunately, these formulations come with non-uniqueness and singularity problems. A widespread strategy is to keep the orientation angles moderate by using updated Lagrangian schemes, \cite{Simo1986, Cardona1988, Ibrahimbegovic1997}. These schemes require rod-specific update procedures and are quite delicate to handle. Crisfield and Jeleni\'c \cite{Crisfield1999} were the first who recognized that most of the updates presented at their time lead to path-dependent solutions.
Alternatively, $\SO(3)$ can also be parameterized as a submanifold, i.e., using more than three coordinates but adding additional conditions. For instance, unit quaternions are given by a quadruple with Euclidean norm one. The 9-parameter method takes all entries of the orthogonal matrix together with the six conditions of orthonormality, \cite{Betsch2002}. Not surprisingly, for most of the parametrizations a rod finite element formulation exists. See \cite{Jelenic1999} for rotation angles, \cite{Gosh2008a} for quaternions, and \cite{Betsch2002, Romero2002} for 9-parameter method, to name a few. 

Another pitfall lies in the approximation of the cross-section orientations. Then not every interpolation strategy of nodal cross-section orientations does preserve the objectivity of the continuous formulation. Assume two different inertial bases that are related by a constant transformation matrix. Then, the nodal cross-section orientations can be expressed with respect to one or the other inertial basis. An interpolation strategy is objective if the interpolated values of these two different nodal cross-section orientations differ only (multiplicatively) by the constant transformation matrix relating the two inertial bases. Before \cite{Crisfield1999}, most interpolation strategies were based on the direct interpolation of the $\SO(3)$-parameters or their incremental updates \cite{Simo1986, Cardona1988, Ibrahimbegovic1995}, which violates objectivity. As a result of non-objective interpolations, the discrete strain measures of a deformed rod change under a superimposed rigid motion. This problem was recognized in \cite{Crisfield1999} and resolved by an interpolation strategy using relative rotation vectors. A few years later, Betsch and Steinmann \cite{Betsch2002} as well as Romero and Armero \cite{Romero2002} proposed simultaneously an objective interpolation strategy that directly interpolates the cross-section base vectors. Admittedly, this strategy comes at the price of abandoning the orthogonality within an element.

The approximations of centerline and cross-section orientations can be considered as the ansatz functions of the rod finite element formulation. The corresponding test functions are the virtual displacement of the centerline and the virtual rotations of the cross-sections. In a Bubnov--Galerkin method, which is the standard in rod finite elements, the test functions follow from a consistent variation of the ansatz functions. In contrast, in a Petrov--Galerkin method \cite{Petrov1940}, the test functions are assumed as independent fields with approximations that are independent of the ansatz functions. Especially, for highly nonlinear ansatz functions, e.g., \cite{Crisfield1999, Jelenic1999}, a Bubnov--Galerkin method becomes almost impenetrable. To simplify the expressions of the discrete virtual work, an independent approximation of the virtual rotation field is suggested in \cite{Jelenic1999}. Very recently, in \cite{Harsch2023a}, the authors of this paper took up the idea and presented a Petrov--Galerkin rod finite element formulation, in which nodal Euclidean transformation matrices are interpolated with the aid of relative twists; a strategy that originates from Sonneville et al.~\cite{Sonneville2014} arising from the $\SE(3)$-structure of the Cosserat rod kinematics. Therefore, this interpolation is called here $\SE(3)$-interpolation.

Another numerical problem in rod theory is known as locking. As discussed in \cite{Balobanov2018}, shear and membrane locking in rods can occur if Kirchhoff (shear-rigidity) and inextensibility constraints follow in the limit case of a parameter tending to zero. This appears for instance for very slender rods if the stiffness parameters are computed in the sense of Saint-Venant by using the material's Young's and shear moduli, respectively, together with the cross-section geometry. Finite elements that are prone to locking cannot fulfill these constraints exactly over the entire element and introduce parasitic dilatation and shear strains. This is either cured by reduced integration, re-interpolation of strain measures \cite{Meier2015, Greco2017}, or mixed formulations \cite{Santos2010, Santos2011, Betsch2016}.

As a consequence of the discussed challenges, a state-of-the-art large-strain rod formulation should be path-independent, objective, and singularity-free. Moreover, a strategy to avoid locking should be available. Often these difficulties are addressed in a very intertwined way. Regularly, rod formulations are proposed together with very specific update rules for static or dynamic numerical analysis \cite{Romero2002}. While the Petrov--Galerkin projection in combination with the $\SE(3)$-interpolation was  already discussed in \cite{Harsch2023a}, the idea of the paper at hand is to present the modularity of this projection method. Indeed, we suggest an entire family of rod finite elements in which the interpolation of the ansatz functions can be exchanged more or less arbitrarily. In particular, the main contributions of the paper are the following:

\begin{itemize}
	\item We present a family of total Lagrangian (thus path-independent) and objective rod finite element formulations parametrized by the nodal total rotation vectors and centerline points. Consequently, the formulations have all a minimal number of six nodal generalized position coordinates.
	\item We apply a Petrov--Galerkin projection method where the virtual centerline displacements are represented in the inertial basis and the virtual rotations are represented in the cross-section-fixed bases. We compare three different well-established interpolation strategies: the $\mR^{12}$-interpolation \cite{Betsch2002, Romero2002}, the $\mR^3 \times \SO(3)$-interpolation \cite{Crisfield1999, Jelenic1999}, and the $\SE(3)$-interpolation \cite{Sonneville2014}.
	\item The introduction of independent velocity fields with the same interpolation as the test functions leads, for the entire family, to the same discrete inertial virtual work functional with a constant and symmetric mass matrix. The discrete equations of motion are in the form of a first-order ordinary differential equation (ODE), which can be integrated using standard methods.
	\item The nodal generalized velocities are coupled with the time-derivative of the nodal generalized position coordinates.
	\item Possible singularities in dynamic simulations are circumvented by employing the concept of the complement rotation vector.
\end{itemize}

Since the interpolation strategies chosen in this paper are from well-established formulations, we want to highlight here some notable differences between these formulations and those presented here. The list is not exhaustive, but it would be too technical for an introduction to go into all the differences. The $\mR^{12}$-interpolation proposed in \cite{Romero2002, Betsch2002} is used together with a Bubnov--Galerkin projection. For that the virtual work functionals are reformulated as functionals of the virtual displacement of the centerline and the variations of the cross-section base vectors, also called directors. While \cite{Betsch2002} introduces additional constraint forces to guarantee orthonormality of the base vectors at the nodes, \cite{Romero2002} proposes a consistent update rule. The two approaches coincide if a null space method \cite{Leyendecker2008a} is applied to \cite{Betsch2002}. The $\mR^3 \times \SO(3)$-interpolation was suggested in combination with a Petrov--Galerkin projection. In contrast to the approach proposed here, in \cite{Jelenic1999} the virtual rotations are represented with respect to the inertial basis. Moreover, the generalized velocities follow from a consistent time derivative of the ansatz functions. These choices result in a second-order ODE with a non-symmetric and configuration dependent mass matrix. In \cite{Sonneville2014}, where the $\SE(3)$-interpolation is applied, the rod theory and its numerical treatment with a Bubnov--Galerkin projection is formulated exclusively in the $\SE(3)$ Lie group setting. Therefore, highly specialized Lie group solvers are needed for the numerical solutions in statics and dynamics. Moreover, the Bubnov--Galerkin projection method requires the computation of the $\SE(3)$-tangent map, its inverse, and its derivatives. Strains, generalized velocities, and generalized virtual displacements are introduced as elements of the Lie algebra $\se(3)$. Consequently, the virtual centerline displacements and the centerline velocities are expressed with respect to the cross-section-fixed basis; this is in contrast to our approach, which uses representations with respect to the inertial bases for these quantities.

The remainder of this paper is organized as follows. Section \ref{sec:Model_of_the_ambient_space}, starts with a digression on vector and affine spaces in conjunction with the introduction of a practical notation containing information about different component representations. In Section~\ref{sec:cosserat_rod_theory}, the Cosserat rod theory is briefly recapitulated in variational form and presented directly as a coordinate representation of \cite{Eugster2020b}. The core of the paper is in Section~\ref{sec:petrov_galerkin_finite_element_formulation}, where the entire family of Petrov--Galerkin finite element formulations is presented. The formulation allows to treat most of the introduced challenges independently. The $\SO(3)$-parameterization by total rotation vectors affects only the nodal coordinates, which are introduced in Section~\ref{sec:nodal_coordinates}. The singularity problem within the kinematic differential equations arising with this choice, and the strategy for avoiding it, are discussed in Section~\ref{sec:kinematic_differential_equations}. The three different objective interpolation strategies are presented in the Sections~\ref{sec:r12_interpolation}-\ref{sec:SE3_interpolation}. The approximation of the test functions and the velocity fields, i.e., the interpolation of the generalized virtual displacements and velocities are treated in Section~\ref{sec:interpolation_of_generalized_virtual}. The family of rod finite element formulations readily follows in Section~\ref{sec:discrete_virtual_work_functionals} by inserting the discrete kinematics into the continuous formulation from Section~\ref{sec:cosserat_rod_theory}. Eventually, in Section~\ref{sec:equations_of_motion_and_static_equilibrium}, the discretization results in the discrete equations of motion of the rod in the form of a first-order ODE. The rod with the $\SE(3)$-interpolation has already been successfully tested against analytical solutions, see \cite{Harsch2023a}. Therefore, Section~\ref{sec:numerical_experiments} gives an in-depth analysis concerning locking and convergence behavior of the different rod formulations. Moreover, a highly dynamic problem with the flexible heavy top is shown, which in the limit of infinitely high stiffness parameters leads to the precession motion of a rigid heavy top. Conclusions are drawn in Section~\ref{sec:conclusions}.

\section{Model of the ambient space and some notational preliminaries}
\label{sec:Model_of_the_ambient_space}

We introduce the three-dimensional Euclidean vector space $\mathbb{E}^3$ as an abstract 3-dimensional real inner-product space. In this paper, only right-handed orthonormal bases are considered. The base vectors of a basis $I$, or $I$-basis, are denoted by $\ve^I_x, \ve^I_y,  \ve^I_z \in \mathbb{E}^3$. The triple ${}_I \va = (a_x^I, a_y^I, a_z^I) \in \mR^3$ contains the components of a vector $\va = a_x^I \ve_x^I + a_y^I \ve_y^I + a_z^I \ve_z^I \in \mathbb{E}^3$ with respect to the $I$-basis. Thus, we carefully distinguish $\mR^3$ from the three-dimensional Euclidean vector space $\mathbb{E}^3$. For computations in components,  triples are treated in the sense of matrix multiplication as $\mR^{3\times1}$-matrices, i.e., as ``column vectors''. For another basis $K$, the same vector $\va=a_x^K \ve_x^K + a_y^K \ve_y^K + a_z^K \ve_z^K  \in \mathbb{E}^3$ has different components collected in ${}_K \va = (a_x^K, a_y^K, a_z^K) \in \mR^3$. Since both $I$- and $K$-basis are right-handed orthonormal bases, the relation between the two representations is given by
\begin{equation}\label{eq:transformation_matrix}
	{}_I \va = \vA_{IK} {}_K \va \, , \quad \text{where } \vA_{IK} \in \SO(3)=\{\vA \in \mR^{3 \times 3}| \vA\T \vA = \mathbf{1}_{3 \times 3} \wedge \det \vA = +1\} \, .
\end{equation}
The transformation matrix $\vA_{IK}$ is an element of the special orthogonal group $\SO(3)$ and transforms the components of a vector with respect to the $K$-basis to its representation in the $I$-basis. Note the suggestive notation, where two adjacent letters cancel each other. 

As a model of the ambient space, we introduce the three-dimensional Euclidean point space $\mathcal{E}^3$, which is an affine space modeled on the Euclidean vector space $\mathbb{E}^3$, see \cite{Crampin1987a} for more details. By definition, for any pair of points $P, Q \in \cE^3$, there exists the unique vector $\vr_{PQ} \in \mathbb{E}^3$ such that $Q = P \hat{+} \vr_{PQ}$, where $\hat{+}:\cE^3 \times \mathbb{E}^3 \to \cE^3$ denotes the affine structure of the affine space. Further, we introduce the concept of a frame, which is the set composed of a point $P \in \cE^3$ together with the orthonormal base vectors of some basis of the vector space $\mathbb{E}^3$. Let the $\cK$-frame be given by the set $\cK = \{P, \ve_x^K, \ve_y^K, \ve_z^K\}$. Then, a point $Q \in \cE^3$ can uniquely be described by its representation with respect to the $\cK$-frame, which is given by the Cartesian coordinates contained in the triple ${}_K\vr_{PQ} \in \mR^3$. Introducing the $\cI$-frame as the set $\{O, \ve_x^I, \ve_y^I, \ve_z^I\}$, the coordinates of the point $Q$ in the $\cI$-frame are given by ${}_I \vr_{OQ}$, which relate to the coordinates in the $\cK$-frame by the affine transformation
\begin{equation}\label{eq:affine_trafo}
	{}_I \vr_{OQ} = {}_I\vr_{OP} + \vA_{IK} {}_K\vr_{PQ} \, .
\end{equation}
Using so-called homogenous coordinates, i.e., extending the triple by an entry that is 1, the affine relation \eqref{eq:affine_trafo} can be written as the linear transformation
\begin{equation}\label{eq:SE3_structure}
	\begin{pmatrix}
		{}_I \vr_{OQ} \\ 1
	\end{pmatrix} =\vH_{\cI \cK} \begin{pmatrix}
		{}_K\vr_{PQ} \\ 1
	\end{pmatrix} \, , \quad \text{with }
	\vH_{\cI \cK} =
	\begin{pmatrix}
		\vA_{IK} & {}_I\vr_{OP} \\
		\mathbf{0}_{1 \times 3} & 1
	\end{pmatrix} \, .
\end{equation}
The matrix $\vH_{\cI \cK}$ is called the Euclidean transformation matrix, which transforms the coordinates of a point in the $\cK$-frame to the coordinates in the $\cI$-frame. The Euclidean transformation matrix $\vH_{\cI \cK}$ is an element of the special Euclidean group $\SE(3)$, which is considered here as a Lie subgroup of the general linear group $\GL(4)$ with the matrix multiplication as group operation. While this paper can be read with a rudimentary knowledge of matrix Lie groups, the interested reader
is referred to \cite[Appendix A]{Harsch2023a} for a deeper understanding of the applied concepts. Direct computation readily verifies that the inverse of $\vH_{\mathcal{I} \mathcal{K}}$ is
\begin{equation}\label{eq:homogenous_transformation}
	\vH_{\mathcal{I}\mathcal{K}}^{-1} = \begin{pmatrix}
		\vA_{IK}\T & - \vA_{IK}\T {}_{I} \vr_{OP} \\
		\mathbf{0}_{1 \times 3} & 1
	\end{pmatrix} \, .
\end{equation}

Both points and base vectors can be considered as functions of time $t$ or some other parameters resulting in parameter-dependent coordinates. The coordinates of a moving point $P = P(t) \in \cE^3$ in the inertial $\cI$-frame are given by the time-dependent triple ${}_I \vr_{OP}={}_I \vr_{OP}(t) \in \mR^3$. 
\section{Cosserat rod theory}\label{sec:cosserat_rod_theory}
\subsection{Centerline and cross-section orientations}
\begin{figure}
	\centering
	\includegraphics{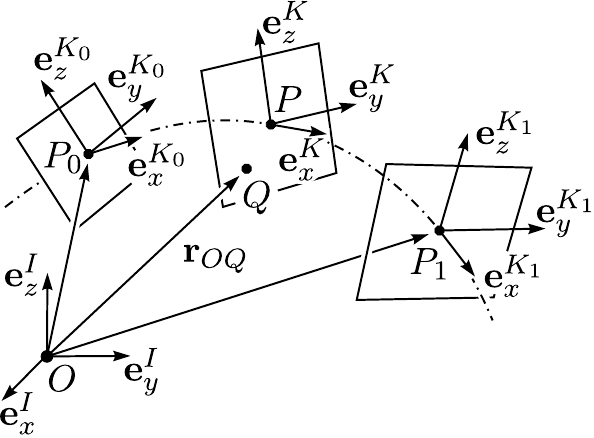}
	\caption{Kinematics of the centerline curve with their attached orthonormal basis vectors.}
	\label{fig:kinematics}
\end{figure}

Let $\xi \in \mathcal{J} = [0, 1] \subset \mR$ denote the centerline parameter and $\mathcal{A}(\xi) \subset \mR^2$ the cross-section area at $\xi$. Considering the rod as a three-dimensional continuum, a point $Q$ of the rod can be addressed by
\begin{equation}\label{eq:rod_kinematics}
	{}_I \vr_{OQ}(\xi,\eta, \zeta, t) = {}_I \vr_{OP}(\xi, t) + \vA_{IK}(\xi, t) \, {}_K \vr_{PQ}(\eta, \zeta) \, , \quad  (\xi ,\ \eta ,\ \zeta) \in \mathcal{B} = \mathcal{J} \times \mathcal{A}(\mathcal{J})\subset \mR^3 \, .
\end{equation}
Herein, ${}_I \vr_{OP}$ are the components with respect to the inertial $I$-basis of the time-dependent centerline curve $\vr_{OP}$, where the subscript $P$ refers to the centerline point, see Figure~\ref{fig:kinematics}. At each centerline point $\vr_{OP}(\xi, t)$, there is a cross-section-fixed $K$-basis determined by the base vectors $\ve_i^K = \ve_i^K(\xi,t),\,i \in \{x,y,z\}$, which are functions of the centerline parameter $\xi$ and time $t$. According to \eqref{eq:transformation_matrix}, the transformation matrix $\vA_{IK}(\xi, t) \in \SO(3)$ relates the representation of a vector in the cross-section-fixed $K$-basis to its representation in the inertial $I$-basis. Consequently, the function $\vA_{IK}$ captures the cross-section orientations, which vary with time and along the rod. Moreover, ${}_K \vr_{PQ}$ denotes the cross-section coordinates, which are independent of time.

The centerline point $P$ together with the cross-section-fixed $K$-basis determine the $\cK$-frame given by the set $\cK = \{P, \ve_x^K, \ve_y^K, \ve_z^K\}$. Comparison of \eqref{eq:rod_kinematics} with \eqref{eq:SE3_structure} reveals the affine structure of the Cosserat rod kinematics and allows to write the motion of the rod \eqref{eq:rod_kinematics} in homogenous coordinates as
\begin{equation}\label{eq:SE3_structure_rod}
	\begin{pmatrix}
		{}_{I} \vr_{OQ} \\
		1
	\end{pmatrix} =
	\vH_{\mathcal{I}\mathcal{K}}
	\begin{pmatrix}
		{}_{K} \vr_{PQ} \\
		1
	\end{pmatrix} \, , \quad \text{with} \quad \vH_{\mathcal{I}\mathcal{K}}(\xi, t) = 
	\begin{pmatrix}
		\vA_{IK}(\xi, t) & {}_I\vr_{OP}(\xi, t) \\
		\mathbf{0}_{1 \times 3} & 1
	\end{pmatrix} \in \SE(3)\, .
\end{equation}
The Euclidean transformation matrix $\vH_{\mathcal{I}\mathcal{K}}(\xi, t)$ relates the coordinates of point $Q$ in the cross-section-fixed $\cK$-frame to the inertial $\cI$-frame. For the application at hand, particularly the group structure of $\vH_{\mathcal{I}\mathcal{K}}(\xi, t)$ will be of relevance. 
\subsection{Velocities, variations, and curvature}
We denote by $\dot{(\bullet)}$ and $(\bullet)_{,\xi}$ the derivatives with respect to time $t$ and centerline parameter $\xi$, respectively. The variation of a function is indicated by $\delta(\bullet)$. The centerline velocity ${}_I \vv_{P}$ and the virtual displacement ${}_I \delta \vr_{P}$ of the centerline are given by the time derivative as well as the variation of the centerline curve
\begin{equation}\label{eq:virtual_displacement_velocity}
	{}_I \vv_{P} = \left({}_I \vr_{OP}\right)^\mathlarger\cdot
	\, , \quad  {}_I \delta \vr_{P} = \delta \left({}_I \vr_{OP}\right) \, .
\end{equation}
The angular velocity of the cross-section-fixed $K$-basis relative to the inertial $I$-basis, in components with respect to the $K$-basis, is defined by
\begin{equation}
	{}_K \vom_{IK} \coloneqq j^{-1}_{\SO(3)} \left( {}_K \widetilde{\vom}_{IK} \right) \, , \quad \text{with} \quad {}_K \widetilde{\vom}_{IK} \coloneqq \vA_{IK}\T \left(\vA_{IK}\right)^\mathlarger\cdot \, ,
\end{equation}
where $j_{\SO(3)} \colon \mR^3 \to \so(3) = \{\vB \in \mR^{3\times3} | \vB\T = -\vB\}$ is the linear and bijective map such that $\widetilde{\vom} \vr = j_{\SO(3)}(\vom) \vr = \vom \times \vr$ for all $\vom, \vr \in \mR^3$. Analogously, we define the virtual rotation as
\begin{equation}\label{eq:virtual_rotation}
	{}_K \delta\vph_{IK} \coloneqq j^{-1}_{\SO(3)} \big( {}_K \delta \widetilde{\vph}_{IK} \big) \, , \quad \text{with} \quad {}_K \delta \widetilde{\vph}_{IK} \coloneqq \vA_{IK}\T \delta \left(\vA_{IK}\right)
\end{equation}
and the  scaled curvature as
\begin{equation}\label{eq:curvature}
	{}_K \bar{\vka}_{IK} \coloneqq j^{-1}_{\SO(3)} \big( {}_K \widetilde{\bar{\vka}}_{IK} \big) \, , \quad \text{with} \quad {}_K \widetilde{\bar{\vka}}_{IK} \coloneqq \vA_{IK}\T \vA_{IK,\xi} \, .
\end{equation}
\subsection{Objective strain measures}
For the reference centerline curve ${}_I \vr_{OP}^0$, the length of the rod's tangent vector is $J = \|{}_I \vr_{OP, \xi}^0\|$. Thus, for a given centerline parameter $\xi$, the reference arc length $s$ is defined by
\begin{equation}
	s(\xi) \coloneqq \int_{0}^{\xi} J(\bar{\xi}) \, \diff[\bar{\xi}] \, .
\end{equation}
Following \cite{Harsch2020a}, the derivative with respect to the reference arc length $s$ of a function $\vf = \vf(\xi,t) \in \mR^3$ can be defined as
\begin{equation}\label{eq:reference_arc_length_derivative}
	\vf_{,s}(\xi,t) \coloneqq \vf_{,\xi}(\xi,t) /J(\xi) \, .
\end{equation}

The objective strain measures of a Cosserat rod, cf. \cite[Section 8.2 and 8.6]{Antman2005}, are
\begin{equation}\label{eq:continuous_strain_measures}
	{}_K \vka_{IK} = {}_K \bar{\vka}_{IK} / J \quad \text{and} \quad {}_K \vga = {}_K \bar{\vga} / J \, , \quad \text{with} \quad {}_K \bar{\vga} \coloneqq (\vA_{IK})\T {}_I \vr_{OP, \xi} \, ,
\end{equation}
which can be gathered in the six-dimensional tuple $\vvep = ({}_K \vga , \ {}_K \vka_{IK}) \in \mR^6$. Therein, the dilatation and shear strains are captured by ${}_K \vga$, while ${}_K \vka_{IK}$ measures torsion and bending.
\subsection{Internal virtual work}
Without loss of generality, we restrict ourselves to hyperelastic material models where the strain energy density with respect to the reference arc length $W = W({}_K \vga, {}_K \vka_{IK}; \xi)$ depends on the strain measures~\eqref{eq:continuous_strain_measures} and possibly explicitly on the centerline parameter $\xi$. By that, the internal virtual work functional is defined as
\begin{equation}\label{eq:internal_virtual_work1}
	\begin{aligned}
		\delta W^\mathrm{int} \coloneqq& -\int_{\mathcal{J}} \delta W J \diff[\xi] = -\int_{\mathcal{J}} \left\{ \delta ({}_K \bar{\vga})\T {}_K \vn + \delta ({}_K \bar{\vka}_{IK})\T {}_K \vm \right\} \diff[\xi] \, ,
	\end{aligned}
\end{equation}
where we have introduced the constitutive equations
\begin{equation}\label{eq:constitutive_equations}
	{}_K \vn \coloneqq \left(\pd{W}{{}_K \vga}\right)\T \, , \quad {}_K \vm \coloneqq \left(\pd{W}{{}_K \vka_{IK}}\right)\T \, .
\end{equation}
Note that even in the inelastic case, where no strain energy density $W$ is available, the internal virtual work~\eqref{eq:internal_virtual_work1} can be used, with internal forces and moments ${}_K \vn$ and ${}_K \vm$ given by different constitutive laws, \cite{Eugster2020b}. Evaluation of $\delta ({}_K \bar\vga)$ and $\delta ({}_K \bar\vka_{IK})$ (see \cite[Appendix B]{Harsch2023a} for details) within the internal virtual work functional~\eqref{eq:internal_virtual_work1} leads to the following expression
\begin{equation}\label{eq:internal_virtual_work2}
	\delta W^\mathrm{int} = -\int_{\mathcal{J}} \big\{ 
	({}_I \delta \vr_{P,\xi})\T \vA_{IK} {}_K \vn + \big[({}_K \delta \vph_{IK})_{,\xi}\big]\T  {}_K \vm 
	- ({}_K \delta\vph_{IK})\T \left[ {}_K \bar{\vga} \times {}_K \vn + {}_K \bar{\vka}_{IK} \times {}_K \vm \right] 
	\big\} \diff[\xi] \, .
\end{equation} 
As in \cite[Equation (2.10)]{Simo1986}, we introduce the diagonal elasticity matrices $\vC_{\vga} = \operatorname{diag}(k_\mathrm{e}, k_\mathrm{s}, k_\mathrm{s})$ and $\vC_{\vka} = \operatorname{diag}(k_\mathrm{t}, k_{\mathrm{b}_y}, k_{\mathrm{b}_z})$ with constant coefficients. In the following, the simple quadratic strain energy density
\begin{equation}\label{eq:strain_energy_density}
	W({}_K \vga, {}_K \vka_{IK}; \xi) = \frac{1}{2} \left({}_K \vga - {}_K \vga^0\right)\T \vC_{\vga} \left({}_K \vga - {}_K \vga^0\right) + \frac{1}{2} \left({}_K \vka_{IK} - {}_K \vka_{IK}^0\right)\T \vC_{\vka} \left({}_K \vka_{IK} - {}_K \vka_{IK}^0\right)
\end{equation}
is used, where the superscript $0$ refers to the evaluation in the rod's reference configuration.
\subsection{External virtual work}
Assume the line distributed external forces ${}_I \vb = {}_I \vb(\xi,t) \in \mR^3$ and moments ${}_K \vc ={}_K \vc(\xi,t) \in \mR^3$ to be given as densities with respect to the reference arc length. Moreover, for $i\in\{0,1\}$, point forces ${}_I \vb_i = {}_I \vb_i(t) \in \mR^3$ and point moments ${}_K \vc_i = {}_K \vc_i(t) \in \mR^3$ can be applied to the rod's boundaries at $\xi_0=0$ and $\xi_1=1$. The corresponding external virtual work functional is defined as
\begin{equation}\label{eq:external_virtual_work}
	\delta W^\mathrm{ext} \coloneqq \int_{\mathcal{J}} \left\{ ({}_I\delta\vr_{P})\T {}_I \vb + ({}_K \delta \vph_{IK})\T {}_K \vc \right\} J \diff[\xi]
	+ \sum_{i = 0}^1 \left[ ({}_I\delta\vr_{P})\T {}_I \vb_i + ({}_K \delta \vph_{IK})\T {}_K \vc_i \right]_{\xi_i} \, .
\end{equation}
\subsection{Inertial virtual work}
Let $\rho_0 = \rho_0(\xi)$ denote the rod's scalar-valued mass density per unit reference volume and $\diff[A]$ the cross-section surface element. It is convenient to define the following abbreviations
\begin{equation}\label{eq:A_S_I_rho}
	A_{\rho_0}(\xi) \coloneqq \int_{\mathcal{A}(\xi)} \rho_0 \, \diff[A] \, , \quad 
	{}_K \vI_{\rho_0}(\xi) \coloneqq \int_{\mathcal{A}(\xi)}  {}_K \widetilde{\vr}_{PQ} \, ({}_K \widetilde{\vr}_{PQ})\T \rho_0 \,  \diff[A] \, .
\end{equation}
Further, using the mass differential $\diff[m] = \rho_0 J \diff[A] \diff[\xi]$ and expressing the variation and second time derivative of $\vr_{OQ}$ in terms of the rod's kinematics \eqref{eq:rod_kinematics}, in case $\vr_{OP}$ is the line of centroids, the inertial virtual work functional of the Cosserat rod can be written as (cf. \cite[Equation (9.54)]{Eugster2020b} for a coordinate-free version)
\begin{equation}\label{eq:inertia_virtual_work}
	\delta W^\mathrm{dyn} \coloneqq -\int_{\mathcal{B}} ({}_I\delta \vr_{OQ})\T {}_I\ddot{\vr}_{OQ} \, \diff[m]= -\int_{\mathcal{J}}
	\big\{({}_I \delta \vr_{P} )\T A_{\rho_0} ({}_I\vv_p)^\mathlarger{\cdot} + ({}_K \delta \vph_{IK})\T
	({}_K \vI_{\rho_0} ({}_K \vom_{IK})^\mathlarger{\cdot} + {}_K \widetilde{\vom}_{IK} {}_K \vI_{\rho_0} {}_K \vom_{IK})\big\} J \diff[\xi]\, .
\end{equation}

\section{Petrov--Galerkin finite element formulation}\label{sec:petrov_galerkin_finite_element_formulation}
\subsection{Lagrangian basis functions}
For the discretization, the rod's parameter space $\mathcal{J}$ is divided into $n_\mathrm{el}$ linearly spaced element intervals $\mathcal{J}^e = [\xi^{e}, \xi^{e+1})$ via $\mathcal{J} = \bigcup_{e=0}^{n_\mathrm{el}-1} \mathcal{J}^e$. For a $p$-th order finite element, the closure of each of the intervals $\mathcal{J}^e$ contains $p + 1$ evenly spaced points $\xi^e_i \in \mathrm{cl}(\mathcal{J}^e) = [\xi^{e}, \xi^{e+1}]$ with $i \in \{0, \dots, p\}$ such that $\xi^e_0 = \xi^e < \xi^e_1 < \dots < \xi^e_p = \xi^{e+1}$. Note, for $e \in \{0, \ldots, n_{\mathrm{el}} -2 \}$, the points $\xi^e_p=\xi^{e+1}_0$ denote the same point $\xi^{e+1}$, which is the boundary point of the adjacent element intervals. It is convenient to use both indexations in the following. For a given element interval $\mathcal{J}^e = [\xi^e, \xi^{e+1})$, the $p$-th order Lagrange basis function and derivative of node $i\in \{0,\dots,p\}$ are
\begin{equation}\label{eq:Lagrangian_polynomials}
	N^{p,e}_i(\xi) = \underset{\substack{0 \leq j \leq p \\ j\neq i}}{\prod} \frac{\xi - \xi^e_j}{ \xi^e_i - \xi^e_j} \, , 
	\quad N^{p,e}_{i,\xi}(\xi) = N_i^{p,e}(\xi) \underset{\substack{k=0 \\ k \neq i}}{\sum^{p}} \frac{1}{\xi - \xi^e_k} \, ,
\end{equation}
where $\xi^e_i$, $\xi^e_j$, and $\xi^e_k$ are the points contained in the set $\{\xi^e_0 = \xi^e, \xi^e_1, \dots, \xi^e_p = \xi^{e+1}\}$.
\subsection{Nodal coordinates}\label{sec:nodal_coordinates}
The here discussed interpolation strategies approximate the centerline curve ${}_I \vr_{OP}$ and the cross-section orientations $\vA_{IK}$ by interpolating nodal centerline points ${}_I\vr_{OP^e_i}(t)\in \mR^3$ and nodal transformation matrices $\vA_{IK^e_i}(t)\in\SO(3)$. For each node $ i \in \{0,\dots,p\}$ within element $e \in \{0,\dots, n_\mathrm{el}-1\}$, it will hold that ${}_I\vr_{OP^e_i}(t) = {}_I \vr_{OP}(\xi^e_i,t)$ and $\vA_{IK^e_i}(t) = \vA_{IK}(\xi^e_i,t)$. Since the nodal transformation matrices $\vA_{IK^e_i}$ are elements of $\SO(3)$, which is a three-dimensional submanifold  of the general linear group $\GL(3)$, a parametrization is required. To obtain a minimal number of nodal coordinates, we propose nodal total rotation vectors $\vps^e_i(t)\in \mR^3$, which parameterize the orientations using the Rodrigues' formula
\begin{equation}\label{eq:SO3_exponential_map}
	\vA_{IK^e_i} = \Exp_{ \SO(3)}(\vps^e_i)=\mathbf{1}_{3\times3} + \frac{\sin(\|\vps^e_i\|)}{\|\vps^e_i\|}\widetilde{\vps}^e_i + \frac{1 - \cos(\|\vps^e_i\|)}{\|\vps^e_i\|^2}(\widetilde{\vps}^e_i)^2 \in \SO(3) \, .
\end{equation}
Note the introduction of a ``capitalized'' $\SO(3)$-exponential map $\Exp_{ \SO(3)}:\mR^3 \to \SO(3)$, which is defined by $\Exp_{ \SO(3)} \coloneqq \exp_{\SO(3)} \circ j_{\SO(3)}$, where the lower case counterpart $\exp_{\SO(3)}\colon \so(3) \to \SO(3)$ denotes the actual $\SO(3)$-exponential map. In the following, we will use this abuse in naming for other Lie group mappings as well. Since \eqref{eq:SO3_exponential_map} has a removable singularity at $\vps^e_i = \mathbf{0}$, for small angles $\|\vps^e_i\| < \epsilon$ \footnote{Here and for all subsequent first-order approximations a critical value of $\epsilon = 10^{-6}$ was used.}, it is beneficial to use the first-order approximation $\Exp_{ \SO(3)}(\vps^e_i)=\mathbf{1}_{3\times3} + \widetilde{\vps}^e_i \in \SO(3)$.

Accordingly, the $N=(p n_\mathrm{el} + 1)$ nodal generalized position coordinates $\vq^e_i(t) = ({}_I \vr_{OP^e_i}, \vps^e_i)(t) \in \mR^6$ are given by the nodal centerline points ${}_I \vr_{OP^e_i}$ and the nodal total rotation vectors $\vps^e_i$ resulting in $n_{\vq} = 6N$ positional degrees of freedom of the discretized rod. The nodal quantities can be assembled in the global tuple of generalized position coordinates $\vq(t) = \big(\vq^0_0, \dots, \vq^0_{p-1}, \dots, \vq^e_{0}, \dots, \vq^e_{p-1}, \dots, \vq^{n_\mathrm{el}-1}_0, \dots,\vq^{n_\mathrm{el}-1}_{p-1}, \vq^{n_\mathrm{el}-1}_p\big)(t) \in \mR^{n_{\vq}}$. For $e \in \{0, \ldots, n_{\mathrm{el}} -2 \}$, the coordinates $\vq^e_p=\vq^{e+1}_0$ refer to the same nodal coordinates. Introducing an appropriate Boolean connectivity matrix $\vC_e \in \mR^{6(p+1) \times n_{\vq}}$, the element generalized coordinates $\vq^e(t) = \big(\vq^e_0, \dots, \vq^e_p\big)(t) \in \mR^{6(p+1)}$ can be extracted from the global generalized position coordinates $\vq$ via $\vq^e = \vC_{e} \vq$. With two further Boolean connectivity matrices $\vC_{\vr, i}, \vC_{\vps, i} \in \mR^{3 \times 6(p+1)}$, the centerline points ${}_I \vr_{OP^e_i}$ and total rotation vectors $\vps^{e}_i$ at node $i$ in element $e$ can be extracted from the element generalized coordinates $\vq^e$ via ${}_I \vr_{OP^e_i} = \vC_{\vr, i} \vq^e$ and $\vps^{e}_i = \vC_{\vps, i} \vq^e$. The Boolean connectivity matrices are only used for the mathematical description of these extraction procedures. During a numerical implementation it is advisable to slice arrays instead of multiply them with Boolean matrices.
\subsection{$\mR^{12}$-interpolation}\label{sec:r12_interpolation}
Following \cite{Betsch2002, Romero2002}, both the centerline and the cross-section orientations are approximated by the piecewise interpolation with $p$-th order Lagrangian polynomials \eqref{eq:Lagrangian_polynomials}$_1$, which can be written as
\begin{equation}\label{eq:R12_interpolation}
	{}_I \vr_{OP}(\xi, \vq) = \sum_{e=0}^{n_\mathrm{el}-1} \chi_{\mathcal{J}^e}(\xi) \sum_{i=0}^p 
	N^{p,e}_i(\xi) {}_I \vr_{OP^e_i}
	\, , \quad
	\vA_{IK}(\xi, \vq) = \sum_{e=0}^{n_\mathrm{el}-1} \chi_{\mathcal{J}^e}(\xi)
	\sum_{i=0}^p 
	N^{p,e}_i(\xi) \vA_{IK^e_i}(\vq)\, ,
\end{equation}
where we have used the characteristic function $\chi_{\mathcal{J}^e} \colon \mathcal{J} \to \{0, 1\}$, which is one for $\xi \in \mathcal{J}^e = [\xi^e, \xi^{e+1})$ and zero elsewhere. In contrast to \cite{Betsch2002, Romero2002}, the nodal cross-section orientations $\vA_{IK^e_i}$ are not coordinates of the system, but are parametrized by the corresponding nodal total rotation vector $\vps^e_i$ in agreement with $\vA_{IK^e_i}(\vq)=\Exp_{ \SO(3)}(\vps^e_i)$. Except at the nodes, the polynomial interpolation \eqref{eq:R12_interpolation}$_2$ leads to non-orthogonal matrices $\vA_{IK}$. This inconsistency with the continuous formulation is accepted as a discretization error that diminishes for decreasing element sizes or increasing polynomial degree. Since \eqref{eq:R12_interpolation} interpolates the three nodal position coordinates and the nine entries of the nodal transformation matrix, the interpolation is called $\mR^{12}$-interpolation.

\begin{figure}[b!]
	\centering
	%
	%
	\begin{subfigure}[b]{.33\textwidth}
		\caption{}
		\vspace{-0.5cm}
		$\mR^{12}$, $p=1$ \\[-1cm]
		\centering
		\includegraphics[width=0.9\textwidth, trim={12cm 10cm 12cm 10cm}, clip]{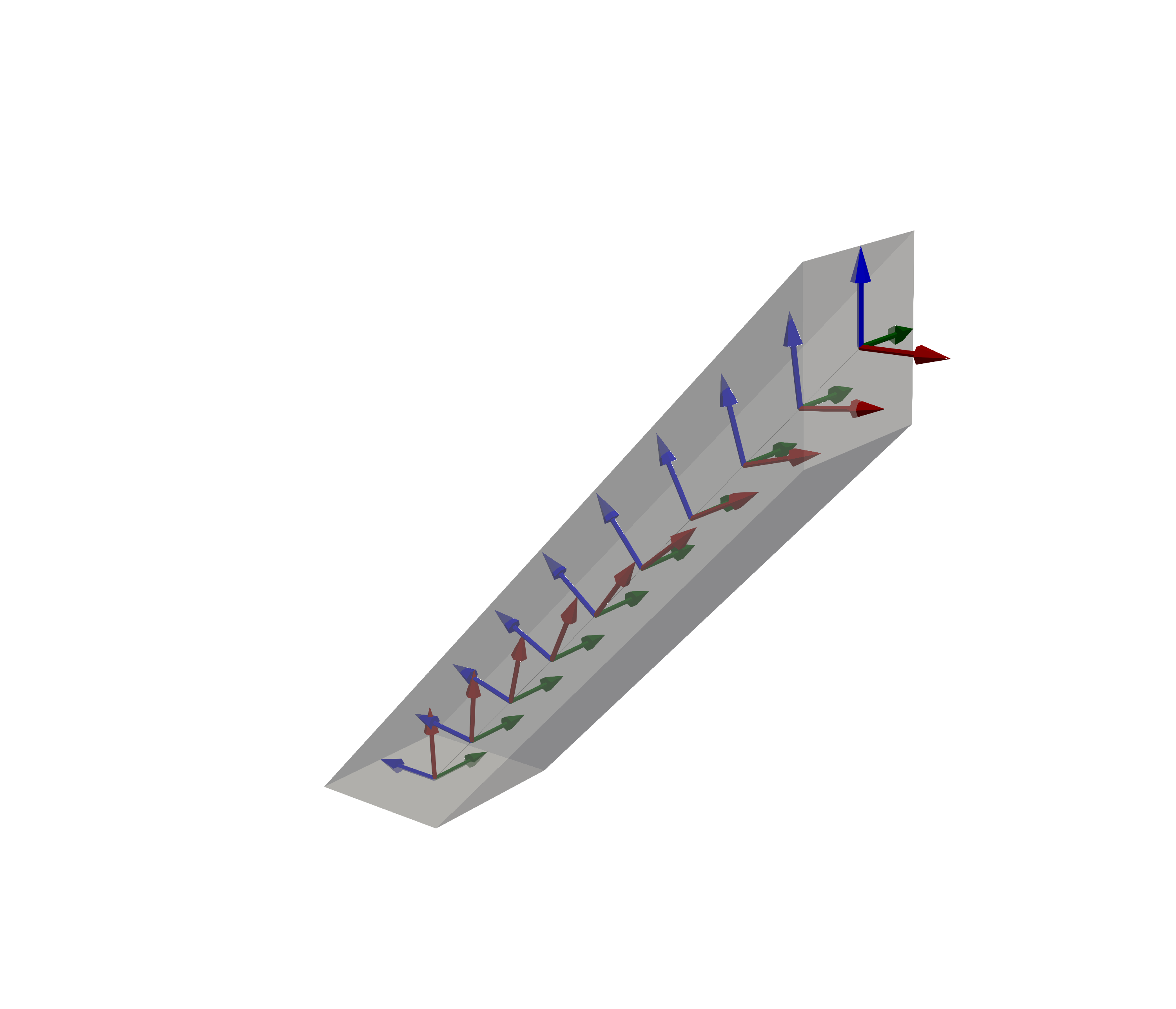}
	\end{subfigure}%
	\begin{subfigure}[b]{.33\textwidth}
		\includegraphics{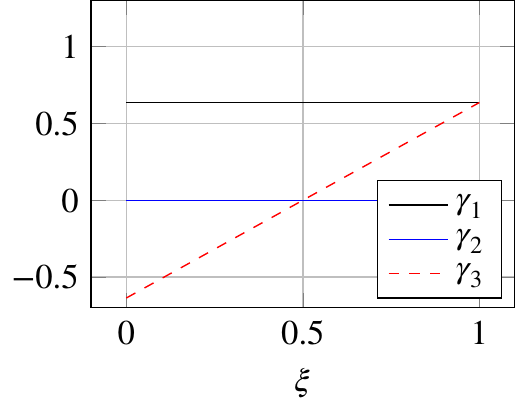}
	\end{subfigure}
	\begin{subfigure}[b]{.33\textwidth}
		\includegraphics{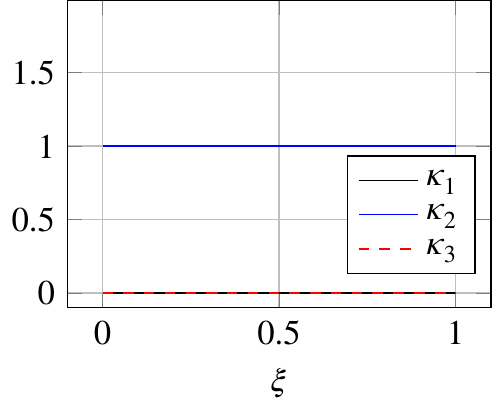}
	\end{subfigure}
	%
	%
	\begin{subfigure}[b]{.33\textwidth}
		\caption{}
		\vspace{-0.5cm}
		$\mR^{12}$, $p=2$ \\[-0.8cm]
		\centering
		\includegraphics[width=\textwidth, trim={12cm 10cm 12cm 10cm}, clip]{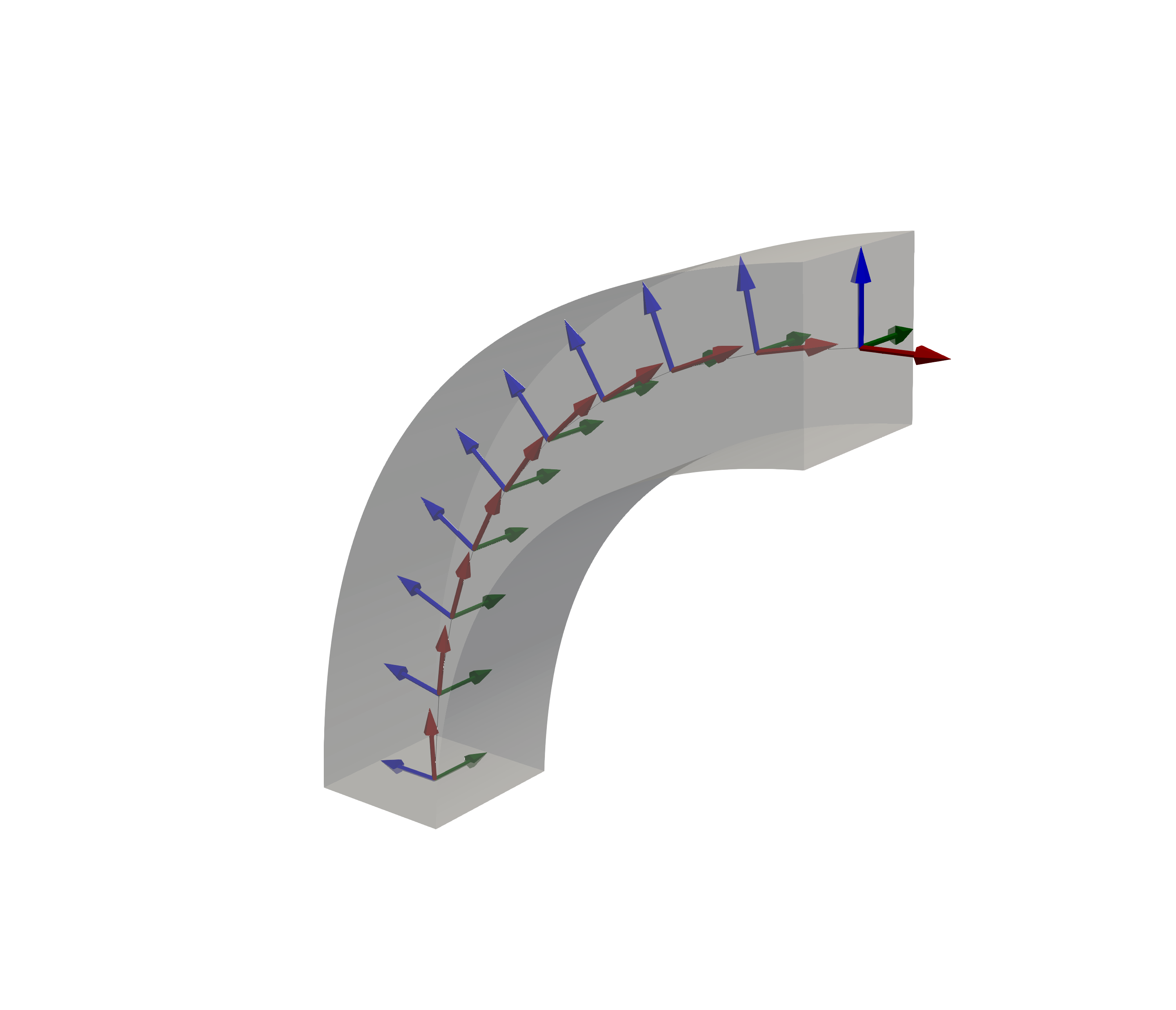}
	\end{subfigure}%
	\begin{subfigure}[b]{.33\textwidth}
		\includegraphics{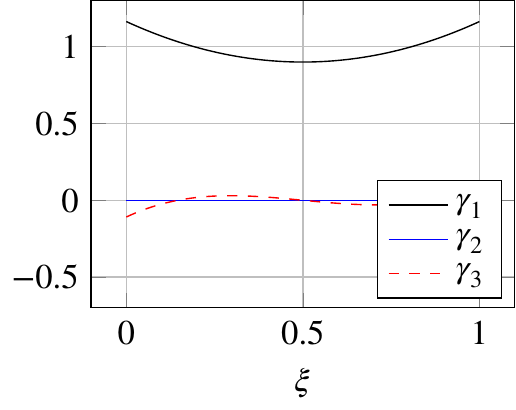}
	\end{subfigure}
	\begin{subfigure}[b]{.33\textwidth}
		\includegraphics{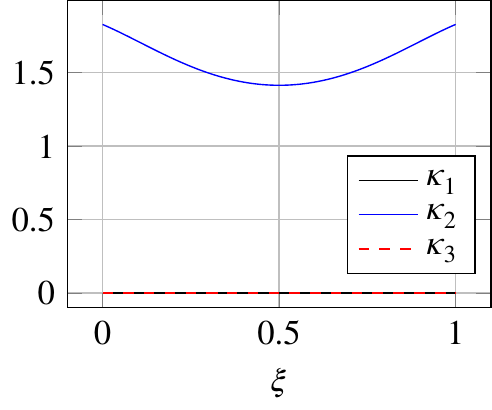}
	\end{subfigure}
	%
	%
	%
	\begin{subfigure}[b]{.33\textwidth}
		\caption{}
		\vspace{-0.5cm}
		$\mR^3 \times \SO(3)$ \\[-1cm]
		\centering
		\includegraphics[width=\textwidth, trim={12cm 10cm 12cm 10cm}, clip]{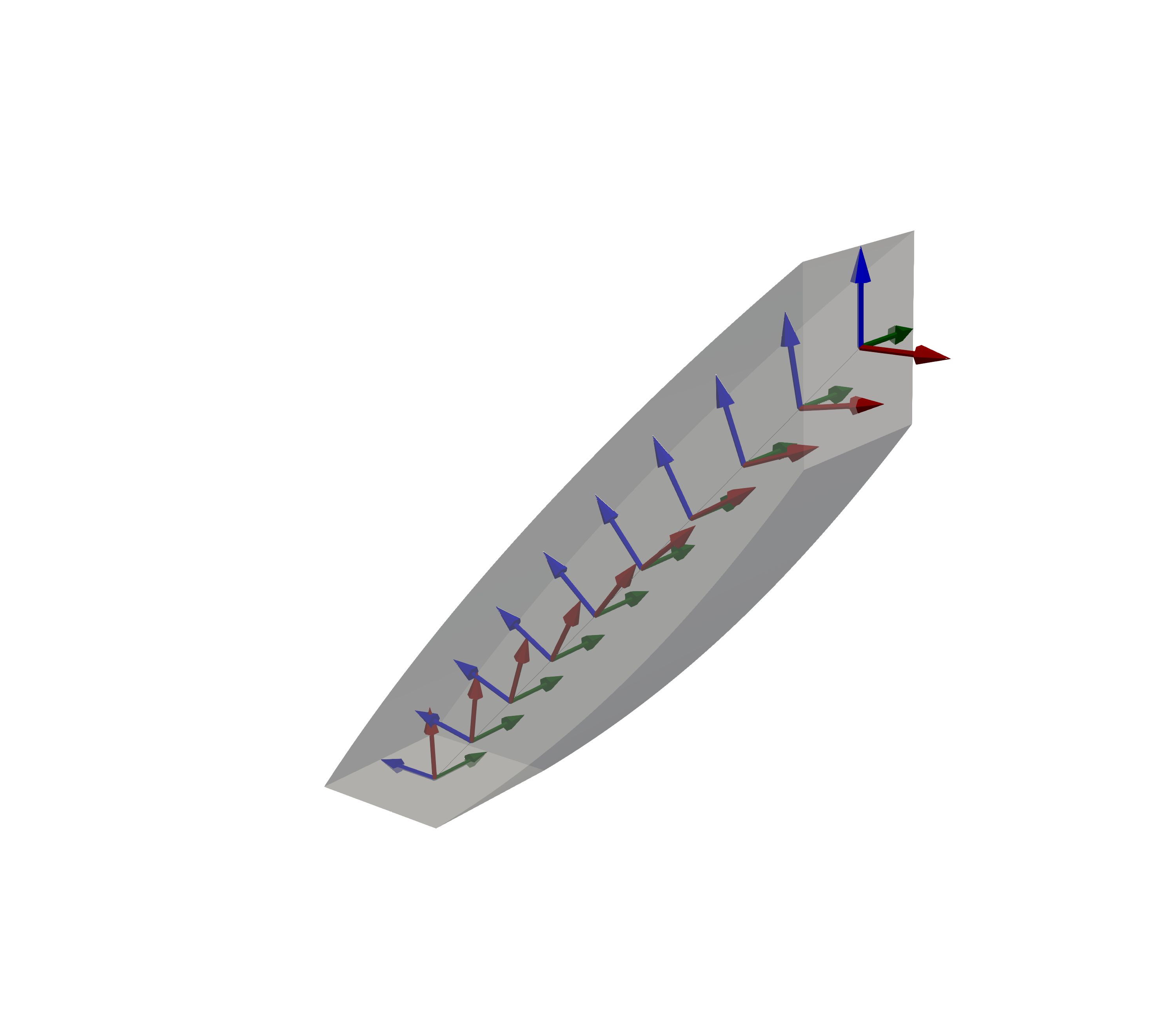}
	\end{subfigure}%
	\begin{subfigure}[b]{.33\textwidth}
		\includegraphics{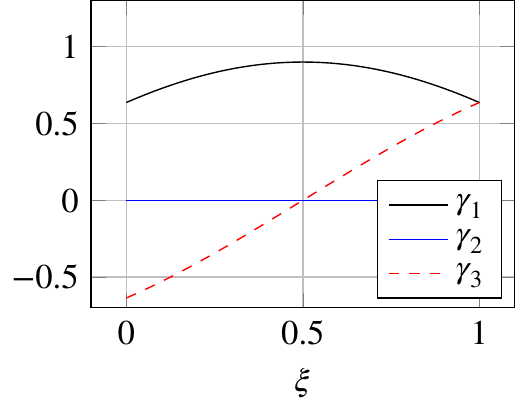}
	\end{subfigure}
	\begin{subfigure}[b]{.33\textwidth}
		\includegraphics{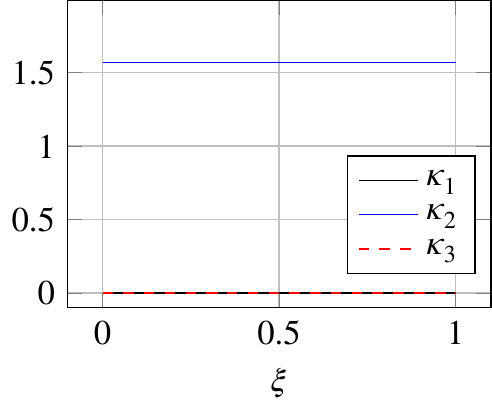}
	\end{subfigure}
	%
	%
	\begin{subfigure}[b]{.33\textwidth}
		\caption{}
		\vspace{-0.45cm}
		$\SE(3)$ \\[-0.8cm]
		\centering
		\includegraphics[width=\textwidth, trim={12cm 10cm 12cm 10cm}, clip]{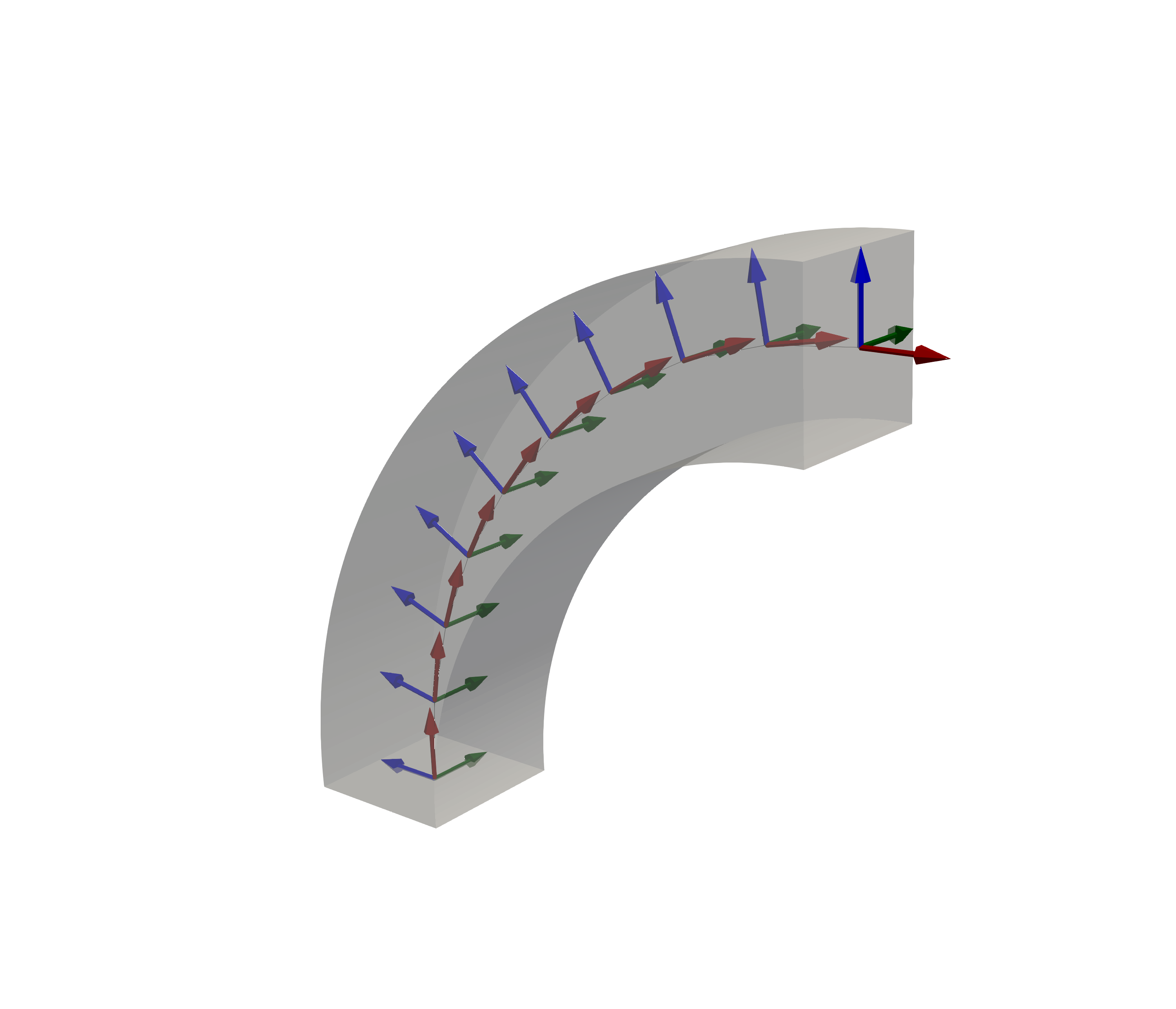}
	\end{subfigure}%
	\begin{subfigure}[b]{.33\textwidth}
		\includegraphics{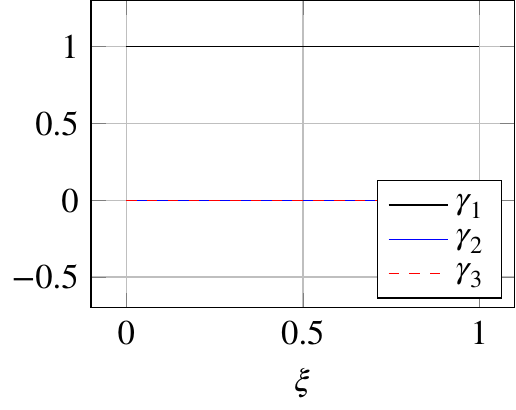}
	\end{subfigure}
	\begin{subfigure}[b]{.33\textwidth}
		\includegraphics{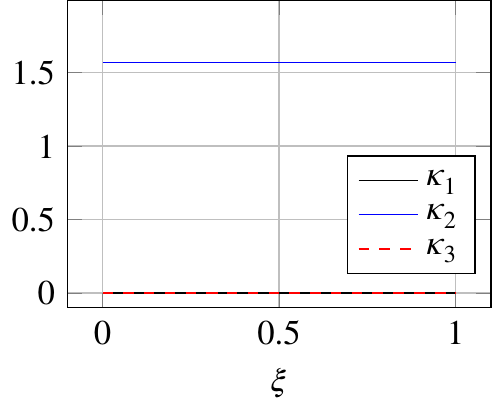}
	\end{subfigure}
	\caption[short]{Volume rendering with base vectors and strain measures of deformed configurations. The nodal rotation vectors are given by $\vps_i = (0, \ \psi_i, \ 0)$ and the nodal centerline by ${}_I \vr_{OP_i} = (1 - \cos(\psi_i), \ 0 , \ \sin(\psi_i)) 2 / \pi$ with $\psi_i = \tfrac{\pi}{2} \tfrac{i}{N - 1}$. \label{fig:strains_quarter_circle}}
\end{figure}

The discretized dilatation and shear strains contained in ${}_K \vga$ are computed by inserting the interpolation \eqref{eq:R12_interpolation} into \eqref{eq:continuous_strain_measures}$_2$. As the interpolated transformation matrices are not orthogonal, the curvature vector cannot be directly extracted by $j_{\SO(3)}^{-1}$. Hence, the discretized version is computed as
\begin{equation}
	{}_K \vka_{IK} = j_{\SO(3)}^{-1}\big(\Skw(\vA_{IK}\T \vA_{IK,\xi})\big)/J \, ,
\end{equation}
where the map $\Skw(\vA) = \frac{1}{2}(\vA - \vA\T) \in \so(3)$ extracts the skew-symmetric part of the matrix $\vA\in\mR^{3\times3}$. See Figure~\ref{fig:strains_quarter_circle}(a) and \ref{fig:strains_quarter_circle}(b) for the discrete strain measures of a quarter circle approximated by one element with $p=1$ and $p=2$, respectively. The objectivity of the discrete strain measures is proven in \cite{Betsch2002}.

%
%
%
\subsection{$\mR^3 \times \SO(3)$-interpolation}\label{sec:r2xSO3_interpolation} 
The objective $\mR^3 \times \SO(3)$-interpolation was originally proposed by Crisfield and Jeleni\'c \cite{Crisfield1999}, who recognized that most rod finite elements at that time do not preserve objectivity after discretization. In our investigations, we will only consider the easiest two-node element of \cite{Crisfield1999}. The centerline is discretized by the piecewise linear interpolation
\begin{equation}\label{eq:R3xSO3_interpolation1}
	{}_I \vr_{OP}(\xi, \vq) = \sum_{e=0}^{n_\mathrm{el}-1} \chi_{\mathcal{J}^e}(\xi) \left(
	N^{1,e}_0(\xi) {}_I \vr_{OP^e_0} + N^{1,e}_1(\xi) {}_I \vr_{OP^e_{1}}
	\right) \, .
\end{equation}
The interpolation of the cross-section orientations is constructed as follows. First, the nodal cross-section orientations $\vA_{IK^e_0}(\vq)= \Exp_{ \SO(3)}(\vps^e_0)$ and $\vA_{IK^e_1}(\vq)= \Exp_{ \SO(3)}(\vps^e_1)$ are evaluated. Second, the relative change of orientation  $\vA_{K^e_0 K^e_{1}} = (\vA_{IK^e_0})\T \vA_{IK^e_{1}}$ is computed, from which the relative rotation vector
\begin{equation}
	\vps^e_{01} = 
	\Log_{\SO(3)}\big(\vA_{K^e_0 K^e_{1}}\big)
\end{equation}
is extracted by using the inverse of $\Exp_{ \SO(3)}$ given by the $\SO(3)$-logarithm map
\begin{equation}\label{eq:Logarithm_SO3}
	\Log_{ \SO(3)}(\vA)= \frac{\omega(\vA)}{2 \sin\big(\omega(\vA)\big)} 
	\begin{pmatrix}A_{32} - A_{23}\\  A_{13} - A_{31}\\  A_{21} - A_{12}\end{pmatrix} \, , \quad \mathrm{with} \quad
	\omega(\vA) = \arccos\left(\frac{1}{2}(\tr(\vA) - 1) \right) \, \text{and } \vA \in \SO(3) \, .
\end{equation}
Note that \eqref{eq:Logarithm_SO3} has a singularity for a relative rotation angle $\omega=\|\vps^e_{01}\|=\pi$. Hence, the interpolation strategy is restricted to applications in which $\|\vps^e_{01}\|<\pi$. A
discretization with a higher number of elements always cures this problem. Lastly, the cross-section orientations are discretized by the ansatz
\begin{equation}\label{eq:R3xSO3_interpolation2}
	\vA_{IK}(\xi, \vq) = \sum_{e=0}^{n_\mathrm{el} - 1} \chi_{\mathcal{J}^e}(\xi) \vA_{IK^e_0}(\vq) \vA_{K^e_0 K}(\xi,\vq) \, , \quad \text{with } \vA_{K^e_0 K}(\xi,\vq) = \Exp_{\SO(3)} \big(N^{1,e}_1(\xi) \,  \vps^e_{01}(\vq)\big) \, .
\end{equation}
Inside element $e$, the interpolation \eqref{eq:R3xSO3_interpolation2} can be understood as a composition of a reference orientation $\vA_{IK_0^e}$ with a relative change of orientation $\vA_{K^e_0 K}$ scaling with $\xi$ such that $\vA_{K^e_0 K}(\xi^e, \vq) = \mathbf{1}_{3\times 3}$ and $\vA_{K^e_0 K}(\xi^{e+1}, \vq) = \vA_{K^e_0 K^e_1}$.

While again ${}_K \vga$ is approximated by inserting the interpolations \eqref{eq:R3xSO3_interpolation1} and \eqref{eq:R3xSO3_interpolation2} into \eqref{eq:continuous_strain_measures}$_2$, in \cite{Crisfield1999} it is shown that the discretized curvature simplifies to the piecewise constant curvature
\begin{equation}\label{eq:discretized_picewise_curvature_measures}
	{}_K \vka_{IK}(\xi, \vq) = \sum_{e=0}^{n_\mathrm{el} - 1} \chi_{\mathcal{J}^e}(\xi) \frac{\vps^e_{01}(\vq)}{\xi^{e+1} - \xi^e} \frac{1}{J}\, .
\end{equation}
Consequently, this interpolation can represent the constant curvatures of the quarter circle with just one element, see Figure~\ref{fig:strains_quarter_circle}(c). 

\subsection{$\SE(3)$-interpolation}\label{sec:SE3_interpolation}
In the $\SE(3)$-interpolation, originally proposed by Sonneville et al. \cite{Sonneville2014} and recently taken up in \cite{Harsch2023a}, the idea from the $\mR^3 \times \SO(3)$-interpolation to use relative rotation vectors is extended to the interpolation of nodal Euclidean transformation matrices with the aid of relative twists. This leads to a strategy coupling the interpolations of centerline points and cross-section orientations. Also here, only a two-node element is presented.

First, the nodal Euclidean transformation matrices
\begin{equation}
	\vH_{\mathcal{I}\mathcal{K}^e_0}(\vq) = 
	\begin{pmatrix}
		\Exp_{ \SO(3)}(\vps^e_0) & {}_I \vr_{OP^e_0} \\
		\mathbf{0}_{3 \times 1} & 1
	\end{pmatrix} \quad \text{and} \quad
	\vH_{\mathcal{I}\mathcal{K}^e_1}(\vq) = 
	\begin{pmatrix}
		\Exp_{ \SO(3)}(\vps^e_1) & {}_I \vr_{OP^e_1} \\
		\mathbf{0}_{3 \times 1} & 1
	\end{pmatrix}
\end{equation} 
are evaluated. Second, the relative Euclidean transformation $\vH_{\mathcal{K}^e_0 \mathcal{K}^e_1} = (\vH_{\mathcal{I}\mathcal{K}^e_0})^{-1} \vH_{\mathcal{I}\mathcal{K}^e_1}$ is computed, from which the relative twist vector
\begin{equation}\label{eq:relative_twist}
	\vth^e_{01} = 
	\Log_{\SE(3)}\big(\vH_{\mathcal{K}^e_0 \mathcal{K}^e_1}\big)=\begin{pmatrix}
		\vT_{\SO(3)}^{-\mathrm{T}}(\vps^e_{01}) \, {}_{K^e_0} \vr_{P^e_0 P^e_{1}} \\
		\vps^e_{01}
	\end{pmatrix} \quad \text{with } \vps^e_{01} = 
	\Log_{\SO(3)}\big(\vA_{K^e_0 K^e_{1}}\big) 
\end{equation}
is extracted by the $\SE(3)$-logarithm map
\begin{equation}\label{eq:SE3_logarithm_map}
	\Log_{ \SE(3)}\left[\begin{pmatrix}\vA & \vr \\ \mathbf{0}_{1 \times 3} & 1 \end{pmatrix}\right] = \begin{pmatrix}
		\vT_{\SO(3)}^{-\mathrm{T}}\big(\Log_{\SO(3)}(\vA)\big) \vr \\
		\Log_{\SO(3)}(\vA)
	\end{pmatrix} \in \mathbb{R}^6\, , \quad \text{where } \vA \in \SO(3), \vr \in \mR^3 \,.
\end{equation}
The $\SE(3)$-logarithm map \eqref{eq:SE3_logarithm_map} requires the inverse of the $\SO(3)$-tangent map given as
\begin{equation}\label{eq:T_SO3_inverse}
	\vT_{\SO(3)}^{-1}(\vps) = \mathbf{1}_{3 \times 3} + \frac{1}{2} \widetilde{\vps} + \bigg(1-\frac{\|\vps\|}{2} \cot\bigg(\frac{\|\vps\|}{2}\bigg)\bigg) \frac{\widetilde{\vps}^2}{\|\vps\|^2} \, .
\end{equation}
As the extraction of the relative twist \eqref{eq:relative_twist} requires the $\SO(3)$-logarithm map \eqref{eq:Logarithm_SO3}, also the $ \SE(3)$-interpolation contains a singularity for a relative rotation angle $\omega=\|\vps^e_{01}\|=\pi$ restricting one element not to be bent more than $180^\circ$. Moreover, the removable singularity of the inverse tangent map at $\vps=\mathbf{0}$ is avoided for $\|\vps\|<\epsilon$ by using the first-order approximation $\vT_{\SO(3)}^{-1}(\vps) = \mathbf{1}_{3 \times 3} + \frac{1}{2}\widetilde{\vps}$.

For the interpolation, we further require the $\SE(3)$-exponential map
\begin{equation}
	\Exp_{\SE(3)} \left[\begin{pmatrix} \vd \\ \vps \end{pmatrix}\right]
	= 
	\begin{pmatrix}
		\Exp_{\SO(3)}(\vps) & \quad \vT_{\SO(3)}\T(\vps) \vd \\
		\mathbf{0}_{1 \times 3} & 1
	\end{pmatrix}\in \SE(3)\, , \text{where } \vd, \vps \in \mR^3 \, ,
\end{equation}
which is determined by the $\SO(3)$-exponential map \eqref{eq:SO3_exponential_map} and the $\SO(3)$-tangent map, see \cite[Section 4.2]{Geradin2001}, 
\begin{equation}
	\vT_{\SO(3)}(\vps) = \mathbf{1}_{3  \times 3} + \left(\frac{ \cos(\|\vps\|) - 1}{\|\vps\|^2}\right) \widetilde{\vps} + \left(\frac{1 - \sin(\|\vps\|)}{\|\vps\|}\right) \frac{\widetilde{\vps}^2}{\|\vps\|^2} \, ,
\end{equation}
with the first-order approximation $\vT_{\SO(3)}(\vps) = \mathbf{1}_{3 \times 3} - \frac{1}{2}\widetilde{\vps}$ for $\|\vps\|<\epsilon$. The interpolation of the nodal Euclidean transformations is then defined as
\begin{equation}\label{eq:relative_interpolation_H_two_node}
	\begin{aligned}
		\vH_{\mathcal{I}\mathcal{K}}(\xi, \vq) &= \sum_{e=0}^{n_\mathrm{el} - 1} \chi_{\mathcal{J}^e}(\xi) \vH_{\mathcal{I}\mathcal{K}^e_0}(\vq)
		\vH_{\cK^e_0 \cK}(\xi, \vq) \, , \quad \text{with } \vH_{\cK^e_0 \cK}(\xi, \vq) = \Exp_{\SE(3)} \big(N^{1,e}_1(\xi) \,  \vth^e_{01}(\vq)\big) \, .
	\end{aligned}
\end{equation}
Inside element $e$, the interpolation \eqref{eq:relative_interpolation_H_two_node} can be understood as a composition of a reference Euclidean transformation $\vH_{\cI \cK_0^e}$ with a relative Euclidean transformation $\vH_{\cK^e_0 \cK}$ scaling with $\xi$ such that $\vH_{\cK^e_0 \cK}(\xi^e, \vq) = \mathbf{1}_{4\times 4}$ and $\vH_{\cK^e_0 \cK}(\xi^{e+1}, \vq) = \vH_{\cK^e_0 \cK^e_1}$. Note the striking resemblance of \eqref{eq:R3xSO3_interpolation2} with \eqref{eq:relative_interpolation_H_two_node}, where just $\SO(3)$- and $\SE(3)$-objects switch places.

According to \eqref{eq:SE3_structure_rod}, the Euclidean transformation matrices $\vH_{\mathcal{I} \mathcal{K}}$ encode the cross-section orientations and the centerline curve. Explicit computation of the components in \eqref{eq:relative_interpolation_H_two_node} leads to the rod's centerline discretization
\begin{equation}\label{eq:SE3_interpolation1}
	{}_I \vr_{OP}(\xi, \vq) = \sum_{e=0}^{n_\mathrm{el} - 1} \chi_{\mathcal{J}^e}(\xi) \big({}_I \vr_{OP^e_0} + N^{1,e}_1(\xi) \, \vA_{IK^e_0}(\vq) \vT_{\SO(3)}\T \big(N^{1,e}_1(\xi) \vps^e_{01}(\vq) \big) \vT_{\SO(3)}^{-\mathrm{T}} \big( \vps^e_{01}(\vq) \big) {}_{K^e_0} \vr_{P^e_0 P^e_{1}}(\vq) \big) 
\end{equation}
and the cross-section orientations
\begin{equation}\label{eq:SE3_interpolation2}
	\vA_{IK}(\xi, \vq) = \sum_{e=0}^{n_\mathrm{el} - 1} \chi_{\mathcal{J}^e}(\xi) \vA_{IK^e_i}(\vq) \Exp_{\SO(3)} \left(N^{1,e}_1(\xi) \,  \vps^e_{01}(\vq)\right) \, ,
\end{equation}
which are interpolated in the same way as in the $\mR^3 \times \SO(3)$-interpolation strategy \eqref{eq:R3xSO3_interpolation2}.

Besides the preservation of objectivity, in \cite{Harsch2023a}, it is shown that the $\SE(3)$-interpolation leads to piecewise constant strains 
\begin{equation}\label{eq:discretized_picewise_strain_measures}
	\vvep(\xi, \vq) = 	\begin{pmatrix}
		{}_K \vga(\xi, \vq) \\
		{}_K \vka_{IK}(\xi, \vq)
	\end{pmatrix} = \sum_{e=0}^{n_\mathrm{el} - 1} \chi_{\mathcal{J}^e}(\xi) \frac{\vth^e_{01}(\vq)}{\xi^{e+1} - \xi^e} \frac{1}{J} \, .
\end{equation}
A fact that has already been recognized in \cite{Sonneville2014}. See Figure~\ref{fig:strains_quarter_circle}(d) for the quarter circle example. Since the piecewise two-node $\SE(3)$-interpolation can exactly represent constant strains within each element, neither membrane nor shear locking will appear with this discretization, see \cite{Harsch2023a} for more details concerning the definition of locking. As we need no further numerical strategies to avoid locking as for instance re-interpolation of strain measures \cite{Meier2015, Greco2017} or mixed formulations \cite{Santos2010, Santos2011, Betsch2016}, this interpolation strategy is called intrinsically locking-free.
\subsection{Interpolation of generalized virtual displacements and velocities}\label{sec:interpolation_of_generalized_virtual}
If applying a Bubnov--Galerkin projection method, each of the just introduced interpolation strategy would lead to a different approximation of the virtual displacements and virtual rotations by inserting  \eqref{eq:R12_interpolation}, \eqref{eq:R3xSO3_interpolation1} and \eqref{eq:R3xSO3_interpolation2}, or \eqref{eq:SE3_interpolation1} and \eqref{eq:SE3_interpolation2} into the definitions of the virtual displacements \eqref{eq:virtual_displacement_velocity}$_2$ and virtual rotations \eqref{eq:virtual_rotation}. Especially for the $\mR^3 \times \SO(3)$- and the $\SE(3)$-interpolations lengthy and cumbersome expressions can be expected in the discretization of the virtual work functionals. Expressions become even longer in a subsequent linearization required for a gradient-based solution strategy such as the common Newton--Raphson method. Therefore, already Jeleni\'c and Crisfield \cite{Jelenic1999} suggested the use of virtual rotation fields that do not follow from a variation of the ansatz function, i.e., the interpolation of the cross-section orientations \eqref{eq:R3xSO3_interpolation2}. Here, we take up again this idea, but formulate the virtual work functionals in terms of virtual rotations expressed in the cross-section-fixed $K$-basis and not in the inertial $I$-basis. 

At the same $N$ nodes as for the nodal generalized position coordinates, we introduce the nodal generalized virtual displacements $\delta \vs^e_i(t) = ({}_I \delta \vr_{P^e_i}, {}_{K^e_i} \delta \vph_{IK^e_i})(t) \in \mR^6$ given by the nodal centerline displacement ${}_I \delta \vr_{P^e_i}(t) \in \mR^3$ and the nodal virtual rotation ${}_{K^e_i} \delta \vph_{IK^e_i}(t) \in \mR^3$. In analogy to the virtual displacements, we also introduce nodal generalized velocities $\vu^e_i(t) = ({}_I \vv^e_{P_i}(t), {}_{K^e_i} \vom_{IK^e_i}(t)) \in \mR^6$ given by the nodal centerline velocity ${}_I \vv_{P^e_i}(t) \in \mR^3$ and the nodal angular velocity ${}_{K^e_i} \vom_{IK^e_i}(t) \in \mR^3$. Similar to the generalized position coordinates $\vq$, the nodal generalized virtual displacements and the nodal generalized velocities are assembled in the global tuple of generalized virtual displacements $\delta \vs(t) \in \mR^{n_{\vq}}$ and global tuple of generalized velocities $\vu(t) \in \mR^{n_{\vq}}$. Again, the Boolean connectivity matrix $\vC_e$ extracts the element virtual displacements $\delta \vs^e(t) = (\delta \vs^e_0, \dots, \delta \vs^e_p)(t) \in \mR^{6(p+1)}$ and the element generalized velocities $\vu^e(t) = (\vu^e_0, \dots, \vu^e_p)(t) \in \mR^{6(p+1)}$ from the global quantities via $\delta \vs^e = \vC_e \delta \vs$ and $\vu^e = \vC_e \vu$. Further, the nodal virtual centerline displacements ${}_I \delta \vr_{P^e_i}$ and centerline velocities ${}_I \vv_{P^e_i}$ can be extracted from the element generalized virtual displacements $\delta \vs^e$ and velocities $\vu^e$ via ${}_I \delta \vr_{P^e_i} = \vC_{\vr, i} \delta \vs^e$ and ${}_I \vv_{P^e_i} = \vC_{\vr, i} \vu^e$, respectively. Identical extraction operations hold for the nodal virtual rotations ${}_{K^e_i} \delta \vph_{IK^e_i} = \vC_{\vps,i} \delta \vs^e$ and angular velocities ${}_{K^e_i} \vom_{IK^e_i}= \vC_{\vps,i} \vu^e$.

In the sense of a Petrov--Galerkin projection~\cite{Petrov1940}, independently of the chosen interpolation strategy for the centerline points and the cross-section orientations, the nodal virtual displacements and rotations are interpolated by $p$-th order Lagrangian basis functions \eqref{eq:Lagrangian_polynomials} in agreement with
\begin{equation}\label{eq:virtual_displacement_interpol}
	{}_I \delta \vr_{P}(\xi, \delta \vs) = \sum_{e=0}^{n_\mathrm{el} - 1} \chi_{\mathcal{J}^e}(\xi) \sum_{i=0}^p N^{p,e}_i(\xi) {}_I \delta \vr_{P^e_i} \, , \quad
	{}_K \delta \vph_{IK}(\xi, \delta \vs) = \sum_{e=0}^{n_\mathrm{el} - 1} \chi_{\mathcal{J}^e}(\xi)\sum_{i=0}^p N^{p,e}_i(\xi) {}_{K^e_i} \delta \vph_{IK^e_i} \, .
\end{equation}
In order to obtain a constant and symmetric mass matrix in the discretized formulation, see~\eqref{eq:discretized_inertial_virtual_work} below, the velocities are considered as independent fields and are interpolated with the same interpolation as the virtual displacements and rotations. Explicitly, they are interpolated by $p$-th order Lagrangian polynomials
\begin{equation}\label{eq:velocity_interpolation}
	{}_I \vv_{P}(\xi, \vu) = \sum_{e=0}^{n_\mathrm{el} - 1} \chi_{\mathcal{J}^e}(\xi)\sum_{i=0}^p N^{p,e}_i(\xi) {}_I \vv_{P^e_i} \, ,  \quad
	{}_K \vom_{IK}(\xi, \vu) = \sum_{e=0}^{n_\mathrm{el} - 1} \chi_{\mathcal{J}^e}(\xi) \sum_{i=0}^p N^{p,e}_i(\xi) {}_{K^e_i} \vom_{IK^e_i} \, .
\end{equation}
\subsection{Kinematic differential equations}\label{sec:kinematic_differential_equations}
The independent introduction of velocity fields \eqref{eq:velocity_interpolation} requires a coupling between position coordinates $\vq$ and velocity coordinates $\vu$. This coupling is satisfied at each node, where the nodal generalized velocities are related to the time derivative of the nodal generalized position coordinates by the nodal kinematic differential equation
\begin{equation}\label{eq:nodal_kinematic_equation}
	\dot{\vq}^e_i =
	\begin{pmatrix}
		{}_I\dot{\vr}_{OP^e_i} \\
		\dot{\vps}^e_i
	\end{pmatrix} =
	\begin{pmatrix}
		\mathbf{1}_{3 \times 3} & \mathbf{0}_{3 \times 3} \\
		\mathbf{0}_{3 \times 3} & \vT^{-1}_{\SO(3)}(\vps^e_i)
	\end{pmatrix}
	\begin{pmatrix}
		{}_I \vv_{P^e_i} \\
		{}_{K^e_i} \vom_{IK^e_i}
	\end{pmatrix} =
	\vB^e_i(\vq^e_i) \vu^e_i \, .
\end{equation}
See \cite[Section 4]{Geradin2001} for a proof of the kinematic differential equation between angular velocity and time derivative of the total rotation vector. Since the inverse of the  $\SO(3)$-tangent map \eqref{eq:T_SO3_inverse} used in~\eqref{eq:nodal_kinematic_equation} exhibits singularities for $\|\vps\| = k 2 \pi$ with $k = 0, 1, 2, \dots$, we
apply the following strategy to avoid them.
As before, the removable singularity at $\vps=\mathbf{0}$ is avoided for $\|\vps\|<\epsilon$ by using the first-order approximation $\vT_{\SO(3)}^{-1}(\vps) = \mathbf{1}_{3 \times 3} + \frac{1}{2}\widetilde{\vps}$. For $k > 0$, the concept of complement rotation vectors \cite{Cardona1988, Ibrahimbegovic1995} is applied. Due to the Petrov--Galerkin projection, it is sufficient to introduce a nodal update that is performed after each successful time step. This update, which corresponds to a change of coordinates for the orientation parametrization, is given by
\begin{equation}\label{eq:change_of_coordinates}
	\vps = 
	\begin{cases}
		\vps \, , & \|\vps\| \leq \pi \, , \\
		\vps^C = \big(1 - 2 \pi / \|\vps\|\big) \vps \, , & \|\vps\| > \pi \, .
	\end{cases}
\end{equation}
It is easy to see that there is no difference whether the nodal transformation matrix $\vA_{IK^e_i}$ is described by the rotation vector $\vps^e_i$ or by its complement $(\vps^e_i)^C = \big(1 - 2 \pi / \|\vps^e_i\| \big) \vps^e_i$, because  $\Exp_{\SO(3)}(\vps^e_i) = \Exp_{\SO(3)}\big((\vps^e_i)^C\big)$. In \cite{Ibrahimbegovic1995}, it is also shown that ${}_{K_i^e} \vom_{IK^e_i} = \vT_{\SO(3)}(\vps^e_i)\dot{\vps}^e_i =  \vT_{\SO(3)}\big((\vps^e_i)^C\big)\big[({\vps}^e_i)^C\big]^\mathlarger{\cdot}$. Hence, neither the nodal kinematic differential equations \eqref{eq:nodal_kinematic_equation} nor the virtual work functionals must be updated upon a change of coordinates \eqref{eq:change_of_coordinates}. For reasonable time steps a both minimal and singularity-free parametrization of $\SO(3)$ is obtained. Note that the proposed strategy of using the complement rotation vector is only required for dynamic simulation, due to the appearing singularity in \eqref{eq:nodal_kinematic_equation}; for static equilibrium problems only the already discussed singularities within one element must be considered. However, these singularities are never a problem, as a higher number of elements always resolves this problem.
\subsection{Discrete virtual work functionals}\label{sec:discrete_virtual_work_functionals}
With the introduced interpolation strategies for ansatz and test functions, we can now discretize the virtual work functionals. In this paper, we will discuss the following formulations, all of which have structurally the same discrete virtual work functionals. For $p \in \{1, 2\}$, the $\mR^{12}$-interpolation with $p$-th order Lagrangian polynomials is combined with $p$-th order approximations of the generalized virtual displacements and velocities. Both the two-node $\mR^3 \times \SO(3)$- and $\SE(3)$-interpolations allow only for approximations of the generalized virtual displacements and velocities with linear Lagrangian polynomials, i.e., $p=1$ for \eqref{eq:virtual_displacement_interpol} and \eqref{eq:velocity_interpolation}.

Inserting \eqref{eq:virtual_displacement_interpol} together with the corresponding approximations for centerline, cross-section orientations and strain measures into \eqref{eq:internal_virtual_work2}, the continuous internal virtual work is approximated by
\begin{equation}
	\begin{aligned}
		&\delta W^\mathrm{int}(\vq; \delta \vs) = \delta \vs\T \vf^{\mathrm{int}}(\vq) \, , \quad \vf^{\mathrm{int}}(\vq) = \sum_{e=0}^{n_\mathrm{el} - 1} \vC_e\T \vf^{\mathrm{int}}_e(\vC_e \vq) \, , \\
		&\vf^{\mathrm{int}}_e(\vq^e) = -\int_{\mathcal{J}^e} \sum_{i=0}^{p}\Big\{ N^{p, e}_{i,\xi} \vC_{\vr, i}\T \vA_{IK} {}_K \vn + N^{p, e}_{i,\xi} \vC_{\vps, i}\T {}_K \vm 
		-N^{p, e}_{i} \vC_{\vps, i}\T \left({}_K \bar{\vga}\times {}_K \vn + {}_K\bar{\vka}_{IK} \times {}_K \vm \right) \Big\} \diff[\xi] \, ,
	\end{aligned}
\end{equation}
where we have introduced the internal forces $\vf^{\mathrm{int}}$ and their element contribution $\vf^{\mathrm{int}}_e$. Above and in the subsequent treatment, we partly suppress the function arguments, which should be clear from the context. Similarly, the external virtual work~\eqref{eq:external_virtual_work} is discretized by
\begin{equation}
	\begin{aligned}
		&\delta W^\mathrm{ext}(\vq; \delta \vs) = \delta \vs\T \vf^{\mathrm{ext}}(\vq) \, , \quad\\ & \vf^{\mathrm{ext}}(\vq) = \sum_{e=0}^{n_\mathrm{el} - 1} \vC_e\T \vf^{\mathrm{ext}}_e(\vC_e \vq) + \vC_{n_\mathrm{el} - 1}\T \left[\vC_{\vr,p}\T {}_I \vb_1 + \vC_{\vps, p}\T {}_K \vc_1 \right]_{\xi=1}
		+ \vC_0\T \left[\vC_{\vr, 0}\T {}_I \vb_0 + \vC_{\vps, 0}\T {}_K \vc_0 \right]_{\xi=0} \, , \\
		&\vf^{\mathrm{ext}}_e(\vq^e) = \int_{\mathcal{J}^e} \sum_{i=0}^{p}\Big\{ N^{p, e}_{i}  \vC_{\vr, i}\T {}_I \vb + N^{p, e}_{i} \vC_{\vps, i}\T {}_K \vc \Big\} J \diff[\xi]  \, ,
	\end{aligned}
\end{equation}
where we have introduced the external forces $\vf^{\mathrm{ext}}$ with their element contributions $\vf^{\mathrm{ext}}_e$. Finally, inserting \eqref{eq:velocity_interpolation} and \eqref{eq:virtual_displacement_interpol} into the inertial virtual work functional~\eqref{eq:inertia_virtual_work} leads to the discrete counterpart
\begin{equation}\label{eq:discretized_inertial_virtual_work}
	\delta W^\mathrm{dyn}(\vu;\delta \vs) = -\delta \vs\T \left\{\vM \dot{\vu} + \vf^{\mathrm{gyr}}(\vu) \right\} \, ,
\end{equation}
where we have made use of the symmetric and constant mass matrix
\begin{equation}
	\vM= \sum_{e=0}^{n_\mathrm{el} - 1} (\vC_e)\T \vM_e \vC_e \, , \quad
	\vM_e = 
	\int_{\mathcal{J}^e} \sum_{i=0}^{p} \sum_{k=0}^{p} N^{p, e}_{i} N^{p, e}_{k} \Big\{
	A_{\rho_0} (\vC_{\vr, i})\T  \vC_{\vr, k} 
	+ (\vC_{\vps, i})\T {}_K \vI_{\rho_0} \vC_{\vps, k}
	\Big\} J \diff[\xi] \, ,
\end{equation}
and the gyroscopic forces
\begin{equation}
	\vf^{\mathrm{gyr}}(\vu) = \sum_{e=0}^{n_\mathrm{el}-1} \vC_e\T \vf_e^{\mathrm{gyr}}(\vC_e \vu) \, , \quad
	\vf^{\mathrm{gyr}}_e(\vu^e) = \int_{\mathcal{J}^e} \sum_{i=0}^{p} N^{p, e}_i \Big\{\vC_{\vps, i}\T {}_K \widetilde{\vom}_{IK} {}_K \vI_{\rho_0} {}_K \vom_{IK} \Big\} J \diff[\xi]  \, .
\end{equation}
Due to the introduction of an independent velocity field \eqref{eq:velocity_interpolation}, the inertial virtual work contributions are for all two-node formulations identical. Moreover, it should be underlined that the mass matrix is constant because the inertial virtual work functional is formulated with angular velocities and virtual rotations expressed in the cross-section-fixed $K$-bases. 

Element integrals of the form $\int_{\mathcal{J}^e} f(\xi) \diff[\xi]$ arising in the discretized internal, external and gyroscopic forces, as well as in the mass matrix, are subsequently computed using an $m$-point Gauss--Legendre quadrature rule. The number of quadrature points depends on the chosen polynomial degree $p$. Full integration requires $m_\mathrm{full} = \lceil(p + 1)^2 / 2\rceil$ quadrature points, where $\lceil x \rceil = \min\{k \in \mathbb{Z} | k \geq x\}$. As it will turn out in Section~\ref{sec:cantilever}, this choice leads to the well known locking behavior. Thus, only the internal virtual work contributions are evaluated using a reduced number of integration points $m_\mathrm{red} = p$.

A discussion about the conservation properties of the just presented semi-discrete formulation can be found in \cite[Appendix D]{Harsch2023a}. While the discrete total energy and linear momentum of the rod is contained in discrete virtual work functionals, there is no algorithmic access to the discrete angular momentum. Consequently, it is not possible to construct a numerical time integration scheme that preserves the discrete total angular momentum. If such a conservation property is crucial for a specific application, the present rod formulation can easily be modified, see also \cite[Appendix D]{Harsch2023a}. Instead of using the nodal virtual rotations ${}_{K_i} \delta \vph_{IK_i}(t) \in \mR^3$ expressed in the cross-section-fixed basis $K_i$, the nodal virtual rotations ${}_{I} \delta \vph_{IK_i}(t) \in \mR^3$ expressed in the inertial basis $I$ are used. Analogously, the nodal angular velocities ${}_{K_i} \vom_{IK_i}(t) \in \mR^3$ are replaced by ${}_{I} \vom_{IK_i}(t) \in \mR^3$. This comes at the price of a configuration dependent mass matrix.

\subsection{Equations of motion and static equilibrium}\label{sec:equations_of_motion_and_static_equilibrium}
The principle of virtual work states that the sum of all virtual work functional has to vanish for arbitrary virtual displacements~\cite[Chapter 8]{dellIsola2020b}, i.e.,
\begin{equation}
	\delta W^\mathrm{tot} = \delta W^\mathrm{int} + \delta W^\mathrm{ext} + \delta W^\mathrm{dyn} \stackrel{!}{=} 0 \quad \forall \delta \vs, \forall t \, .
\end{equation}
Thus, the equations of motion
\begin{equation}
	\dot{\vu} = \vM^{-1} \left(\vf^\mathrm{gyr}(\vu) + \vf^\mathrm{int}(\vq) + \vf^\mathrm{ext}(\vq)\right)
\end{equation}
have to be fulfilled for each instant of time $t$. Further, the nodal generalized velocities $\vu^e_i$ are related to the time derivatives of the nodal generalized position coordinates $\dot{\vq}^e_i$ via the nodal kinematic differential equation~\eqref{eq:nodal_kinematic_equation}, which can be assembled to a global kinematic differential equation of the form $\dot{\vq} = \vB(\vq) \vu$. Depending on the specific application, either the system of ordinary differential equations
\begin{equation}
	\begin{aligned}
		\dot{\vq} &= \vB(\vq) \vu \, , \\
		\dot{\vu} &= \vM^{-1} \left(\vf^\mathrm{gyr}(\vu) + \vf^\mathrm{int}(\vq) + \vf^\mathrm{ext}(\vq)\right) \, ,
	\end{aligned}
\end{equation}
or the nonlinear generalized force equilibrium
\begin{equation}\label{eq:nonlinear_generalized_force_equilibrium}
	\vf^\mathrm{int}(\vq) + \vf^\mathrm{ext}(\vq) = \mathbf{0}
\end{equation}
is obtained. The system of ordinary differential equations can be solved using standard higher-order ODE solvers (e.g. family of explicit \cite{Hairer1993} and implicit \cite{Jay1995, Hairer2002} Runge--Kutta methods or structure-preserving algorithms \cite{Jay1996, Hairer2006}).  In order to apply well-established methods from structural dynamics like the Newmark-$\beta$ \cite{Newmark1959} or the generalized-$\alpha$ method \cite{Chung1993}, a slightly modified update of the generalized coordinates has to be applied, see Equation (37a) and (37b) of \cite{Arnold2016}. Alternatively, a generalized-$\alpha$ formulation for first-order differential equations \cite{Jansen2000} can be used without any modifications. The nonlinear generalized force equilibrium~\eqref{eq:nonlinear_generalized_force_equilibrium} is solved by any root-finding algorithm, e.g., Newton--Raphson, Riks. Note that a system of linear equations with a non-symmetric matrix must be solved in each iteration.
\begin{table}[b]
	\centering
	\caption{Experimental parameters of the \emph{cantilever experiment}.}\label{tab:cantilever}
	\begin{tabular}{lllll}
		\toprule
		$\rho$ & $\SI{e1}{}$ & $\SI{e1}{}$ & $\SI{e3}{}$ & $\SI{e4}{}$ \\
		$atol$ & $\SI{e-8}{}$ & $\SI{e-10}{}$ & $\SI{e-12}{}$ & $\SI{e-14}{}$ \\
		\bottomrule
	\end{tabular}
\end{table}
\section{Numerical experiments}\label{sec:numerical_experiments}
\subsection{Cantilever experiment}\label{sec:cantilever}
\begin{figure}
	\centering
	\begin{subfigure}[b]{0.23\textwidth}
		\caption{}
		\includegraphics[width=\textwidth, trim=20cm 8cm 20cm 8cm, clip]{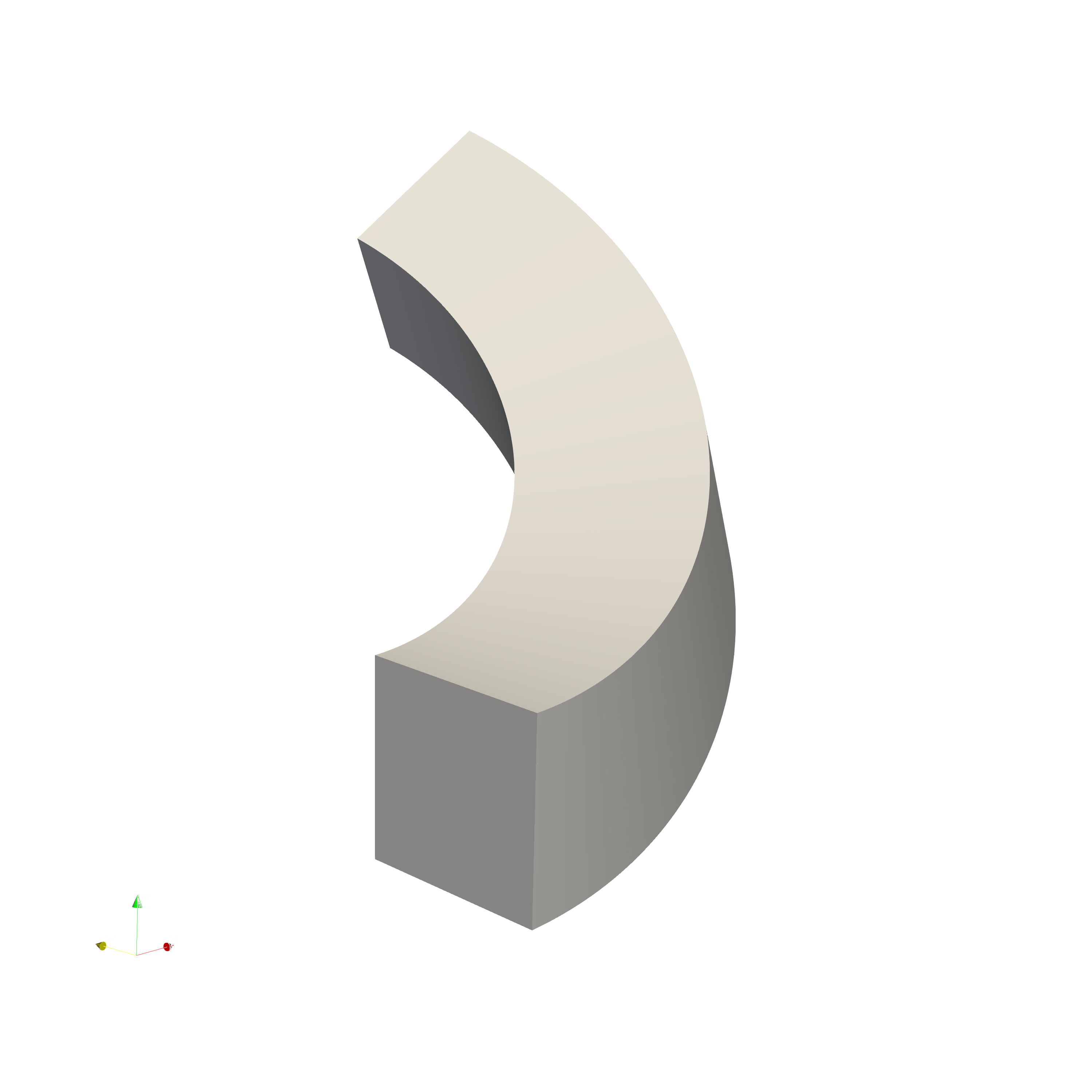}
	\end{subfigure}
	~
	\begin{subfigure}[b]{0.23\textwidth}
		\caption{}
		\includegraphics[width=\textwidth, trim=20cm 8cm 20cm 8cm, clip]{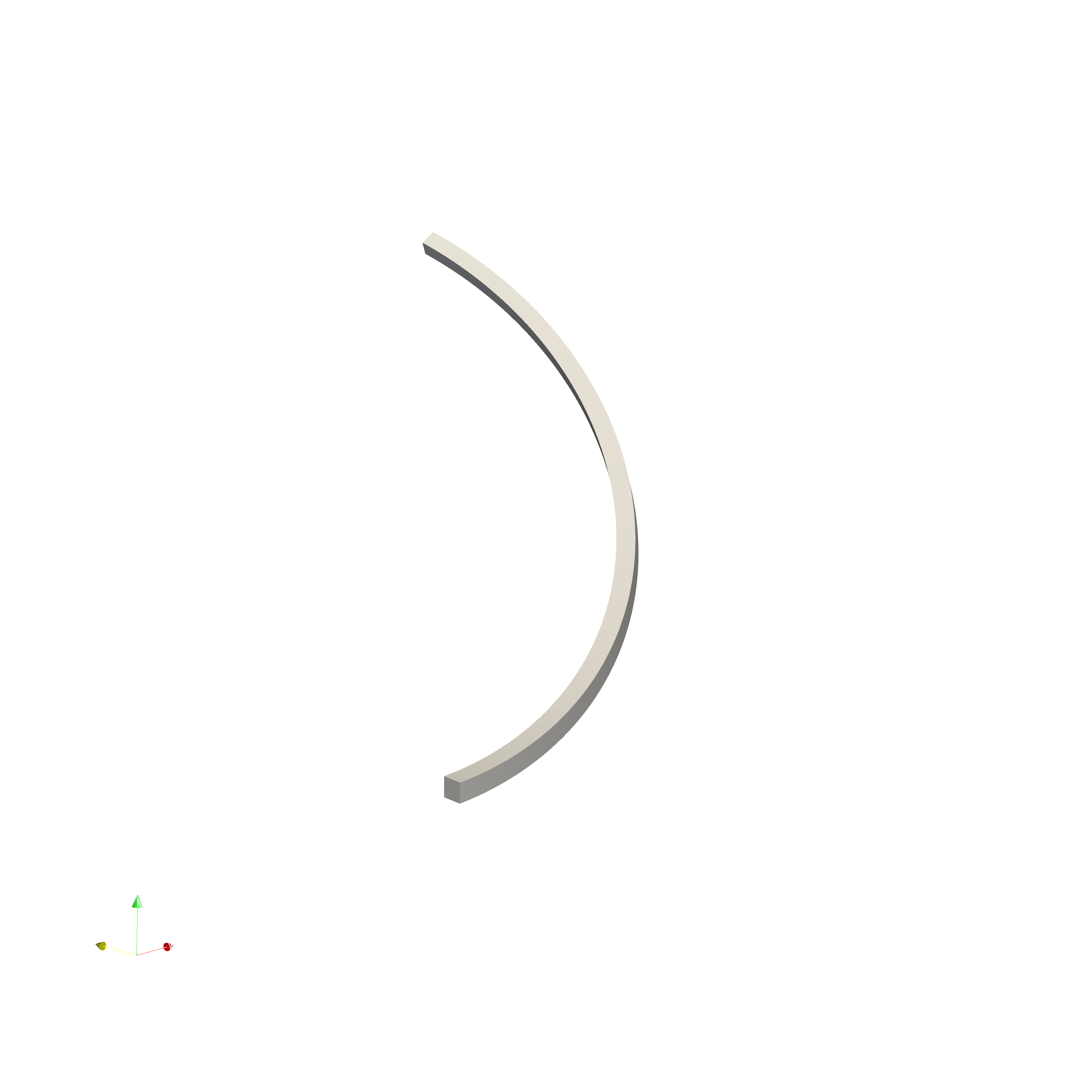}
	\end{subfigure}
	~
	\begin{subfigure}[b]{0.23\textwidth}
		\caption{}
		\includegraphics[width=\textwidth, trim=20cm 8cm 20cm 8cm, clip]{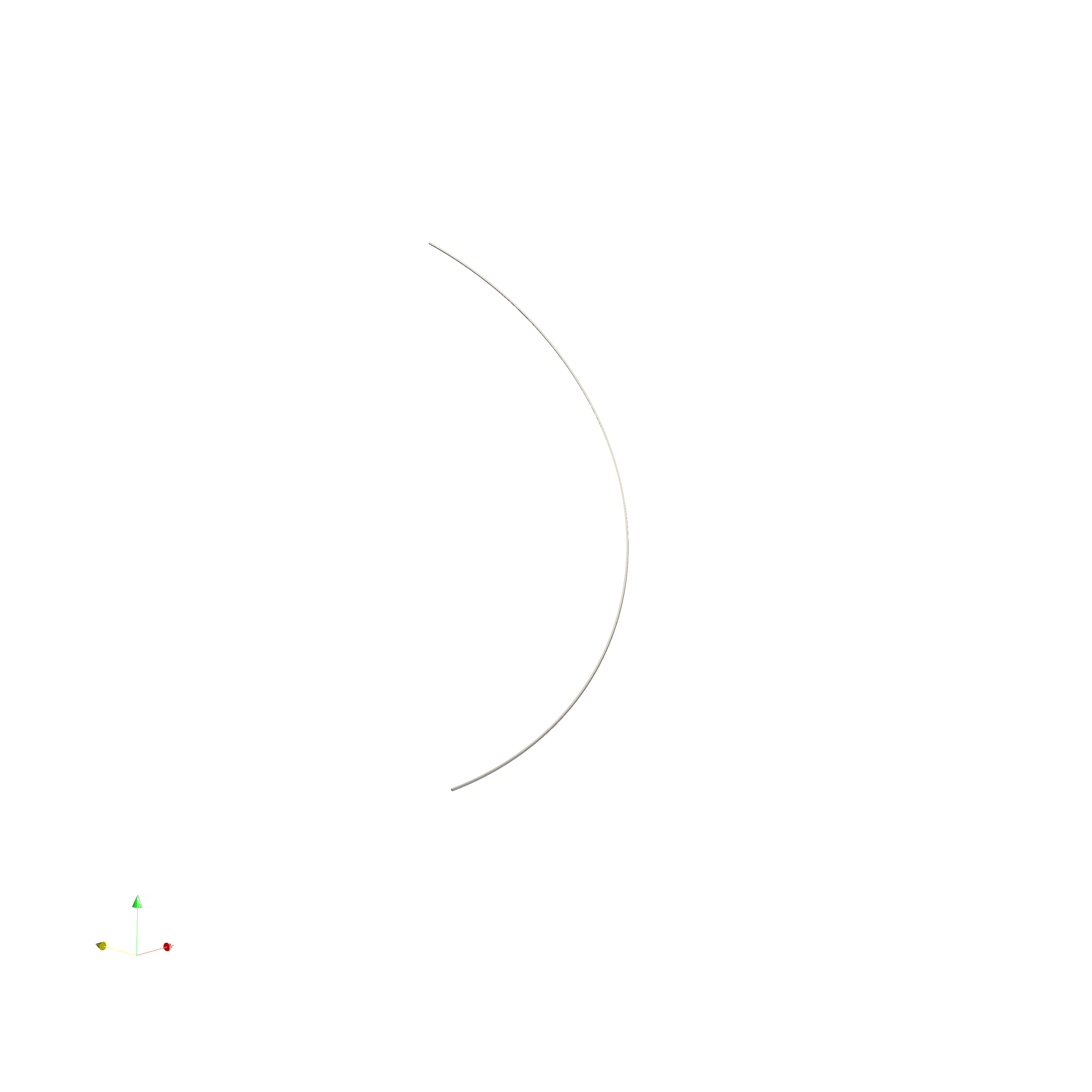}
	\end{subfigure}
	~
	\begin{subfigure}[b]{0.23\textwidth}
		\caption{}
		\includegraphics[width=\textwidth, trim=20cm 8cm 20cm 8cm, clip]{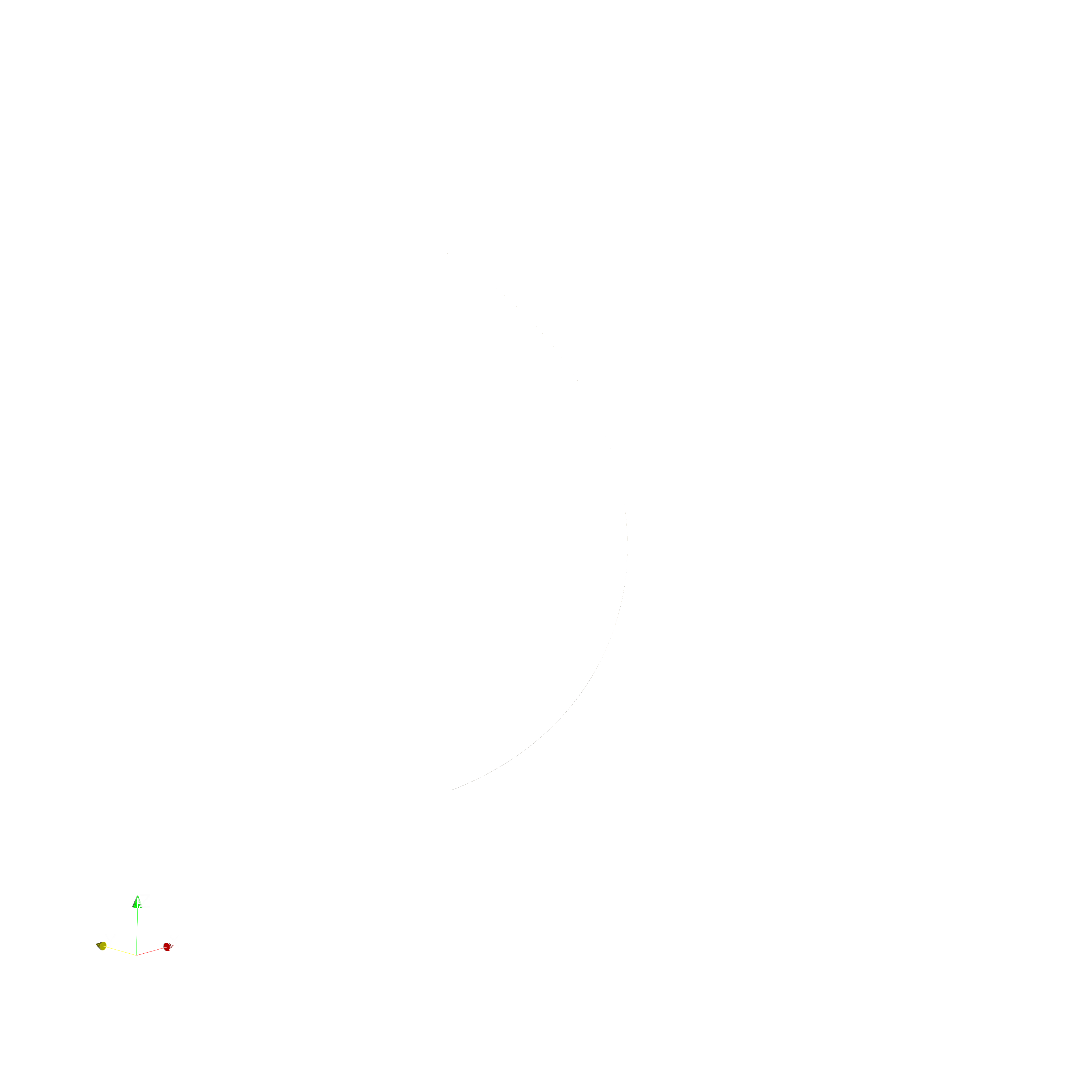}
	\end{subfigure}
	\caption{Visualization of deformed configuration of the \emph{cantilever experiment} for (a) $\rho = \SI{e1}{}$, (b) $\rho = \SI{e2}{}$, (c) $\rho = \SI{e3}{}$ and (d) $\rho = \SI{e4}{}$ for five $\SE(3)$-elements.}
	\label{fig:quarter_circle_final_configurations}
\end{figure}
We consider an initially straight cantilever rod of length $L = \SI{e3}{}$ with a quadratic cross-section of width $w$ subjected to a tip moment ${}_K\vc_1 = (0, \ 0, \ 0.5 \pi k_\mathrm{b} / L)$ and an out-of-plane tip load ${}_I\vb_1 = \vA_{IK} \, (0, \ 0, \ 0.5 \pi k_\mathrm{b} / L^2)$. In order to investigate the presence of locking, different slenderness ratios $\rho = L / w \in \{\SI{e1}{}, \SI{e2}{}, \SI{e3}{}, \SI{e4}{}\}$ are considered, i.e., widths $w \in \{\SI{e2}{}, \ \SI{e1}{}, \ \SI{e0}{}, \ \SI{e-1}{}\}$. Further, the elastic constants are given in terms of the Young's and shear moduli $E=1$ and $G=0.5$. That is, axial stiffness $k_\mathrm{e} = EA$, shear stiffness $k_\mathrm{s} = GA$, bending stiffnesses $k_{\mathrm{b}} = k_{\mathrm{b}_y} = k_{\mathrm{b}_z}= EI$ and torsional stiffness $k_\mathrm{t} = 2GI$, together with $A=w^2$ and $I =w^4/12$. A visualization of the deformed configurations for all different slenderness ratios is shown in Figure~\ref{fig:quarter_circle_final_configurations}.

The required element integrals of the presented rod formulations are integrated using (i) $m_\mathrm{full} = \lceil(p + 1)^2 / 2\rceil$ and (ii) $m_\mathrm{red} = p$ quadrature points. Further, the static equilibrium configurations were found by using a Newton--Raphson method with absolute tolerances $atol$ in terms of the max error of the total generalized forces, given in Table~\ref{tab:cantilever}. The final loads were applied with $50$ load increments. Prescribed kinematic boundary conditions were incorporated into the principle of virtual work using perfect bilateral constraints \cite{Geradin2001}.

As there is no analytical solution for this load case, we require a reference solution whose choice is postponed to the discussion. We introduce the root-mean-square error of relative twists
\begin{equation}
	e^k_{\vth} = \frac{1}{k} \sqrt{\sum_{i=0}^{k - 1} \Delta \vth(\xi_i)\T \Delta \vth(\xi_i) } \, , \quad \Delta \vth(\xi_i) = \Log_{\SE(3)}\big(\vH_{\cI \cK}(\xi_i)^{-1} \vH_{\cI \cK}^{*}(\xi_i)\big) \, , \quad \xi_i = \frac{i}{k - 1}
\end{equation}
as a unified error measure for positions and cross-section orientations.

\begin{figure}
	\centering
	\includegraphics{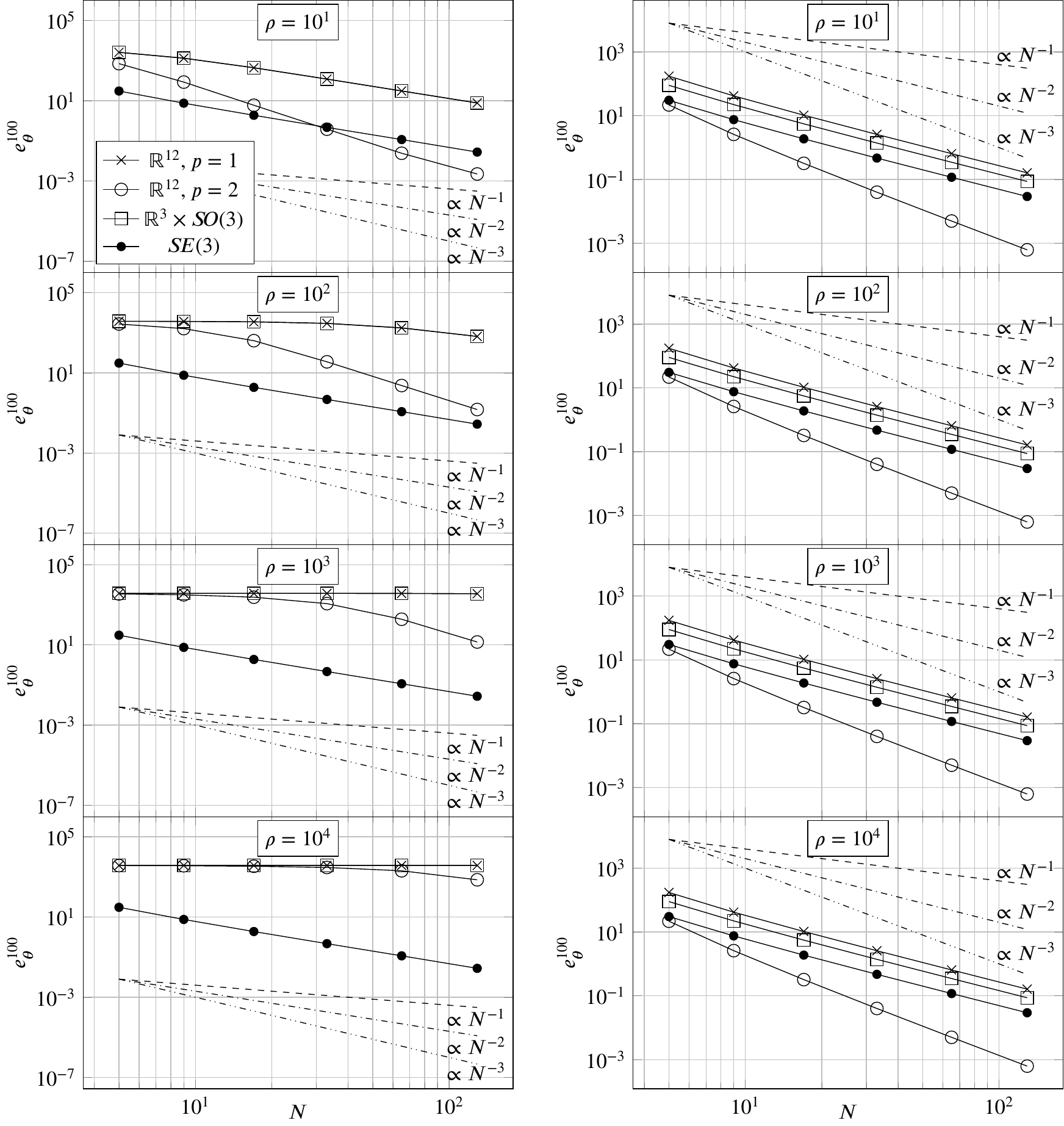}
	\caption{Spatial convergence rates for different interpolation strategies and slenderness ratios using full integration (first column) and reduced integration (second column). }
	\label{fig:convergence_cantilever}
\end{figure}

If no reduced integration is performed, all but the $\SE(3)$-element suffer from locking, see first column of Figure~\ref{fig:convergence_cantilever}. Thus, the solution found by the $\SE(3)$ formulation with $512$ elements ($513$ nodes) was chosen as reference for this experiment. In contrast, the application of reduced integration completely cures the locking phenomenon. Thereby, the second-order $\mR^{12}$-interpolation results in third order spatial convergence, while the other formulations indicate a second-order rate, see second column of Figure~\ref{fig:convergence_cantilever}. Moreover, the quadratic $\mR^{12}$-element outperforms all other formulations by magnitudes using the same number of unknowns. Hence, we use $256$ elements ($513$ nodes) of the quadratic $\mR^{12}$-element as reference for the second experiment.

\begin{figure}
	\begin{subfigure}[b]{.48\textwidth}
		\caption{}
		\centering
		\includegraphics{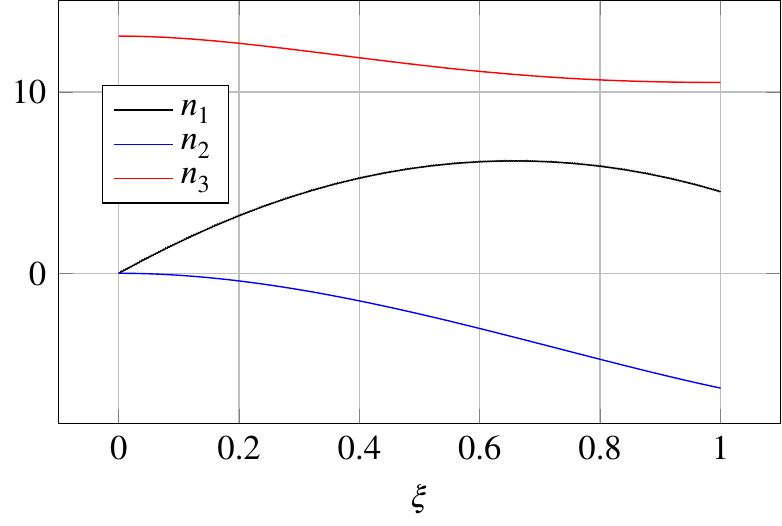}
	\end{subfigure}%
	\begin{subfigure}[b]{.48\textwidth}
		\caption{}
		\centering
		\includegraphics{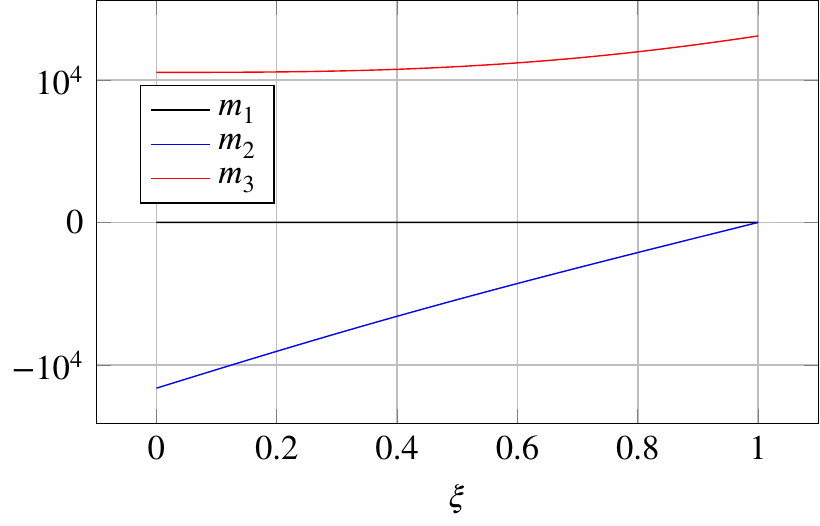}
	\end{subfigure}%
	\caption{Internal forces and moments of the \emph{cantilever experiment} for a moderate slenderness $\rho = \SI{e2}{}$, computed with $256$ elements ($513$ nodes) of the quadratic $\mR^{12}$-element.}
	\label{fig:strain_measures_cantilever}
\end{figure}

For a moderate slenderness $\rho = \SI{e2}{}$, the internal forces and moments of the reference solution are plotted in Figure~\ref{fig:strain_measures_cantilever}. Since the deformation is inhomogeneous, a nonlinear progression is clearly evident.
\subsection{Flexible heavy top}
Inspired by the investigation of Mäkinen \cite{Maekinen2007}, the dynamics of an elastic heavy top is studied here. On the one hand, this example demonstrates the capability of the presented formulations to be solved using standard ODE solvers. On the other hand, in the limit case of an infinite stiff rod, the rod shows the well-known behavior of a heavy top. The motion of the heavy top is described by Euler's equations, see \cite[Equations (1.83) and (3.35)]{Magnus1971}, whose solution is taken as a solution in the subsequent investigation. For high stiffnesses of the rod, we thus expect the solution to be close to the one of a rigid body.

Let the top be given by a cylinder of radius $r=0.1$ and length $L=0.5$ with cross-section area $A = \pi r^2$ and second moment of area $I = \pi r^4 / 4$. The cylinder is subjected to a constant distributed line force ${}_I \vb = \rho_0 A (0, \ 0, \ - 9.81)$. The stiff rod should have a uniform density $\rho_0 = 8000$, Young's modulus $E=\SI{210e6}{}$ and shear modulus $G = E / (2(1 + \nu))$, with a Poisson's ratio $\nu = 1/3$. Consequently, it has an axial stiffness $k_\mathrm{e} = EA$, shear stiffness $k_\mathrm{s} = GA$, bending stiffnesses $k_\mathrm{b} = k_{\mathrm{b}_y} = k_{\mathrm{b}_z} = EI$ and torsional stiffness $k_\mathrm{t} = 2GI$. We also considered a soft rod for which all stiffnesses were scaled by a factor $\SI{2.5e-3}{}$.

Since in the previous example the second-order $\mR^{12}$-interpolation outperforms all other formulations, we restricted the dynamic investigation to this formulation. The top was discretized using a single quadratic $\mR^{12}$-element (3 nodes) and reduced integration was applied for the evaluation of the internal virtual work contributions. The initial position was such that the top points from the origin in positive $\ve_x^I$-direction, i.e., $\vq^0 = (\mathbf{0}_{3\times1}, \ \mathbf{0}_{3\times1}, \ \vr_1, \mathbf{0}_{3\times1}, \ \vr_2, \mathbf{0}_{3\times1})$ with $\vr_1 = (L / 2, \ 0, \ 0)$ and $\vr_2 = (L, \ 0, \ 0)$. Its initial velocities were chosen such that in the case of a rigid rod a perfect precession motion, see \cite[Section 3.3.2 c)]{Magnus1971}, is obtained, i.e., $\vu^0 = (\mathbf{0}_{3\times1}, \ \vOm, \ \vOm \times \vr_1, \vOm, \ \vOm \times \vr_2, \vOm)$ with the angular velocity $\vOm = (\Omega, \ 0, \ \Omega_\mathrm{pr})$, where $\Omega = 50 \pi$ and $\Omega_\mathrm{pr} = g L / (r^2 \Omega)$. Finally, the motion of the top was constrained such that the first node coincides with the origin for all times. This can either be guaranteed by removing the corresponding degrees of freedom from the set of unknowns, or by using the concept of perfect bilateral constraints, \cite{Geradin2001}.

Using a standard fourth-order Runge--Kutta method, with the absolute and relative tolerances $atol = rtol = \SI{1e-8}{}$, see \cite{Hairer1993}, the simulations were performed until a final time of $t_1 = 2 \pi / \Omega_\mathrm{pr}$ was reached, i.e., the rigid top performed a full rotation. 

\begin{figure}[h!]
	\centering
	\begin{subfigure}[b]{0.47\textwidth}
		\caption{}
		\vspace{-0.15cm}
		\includegraphics{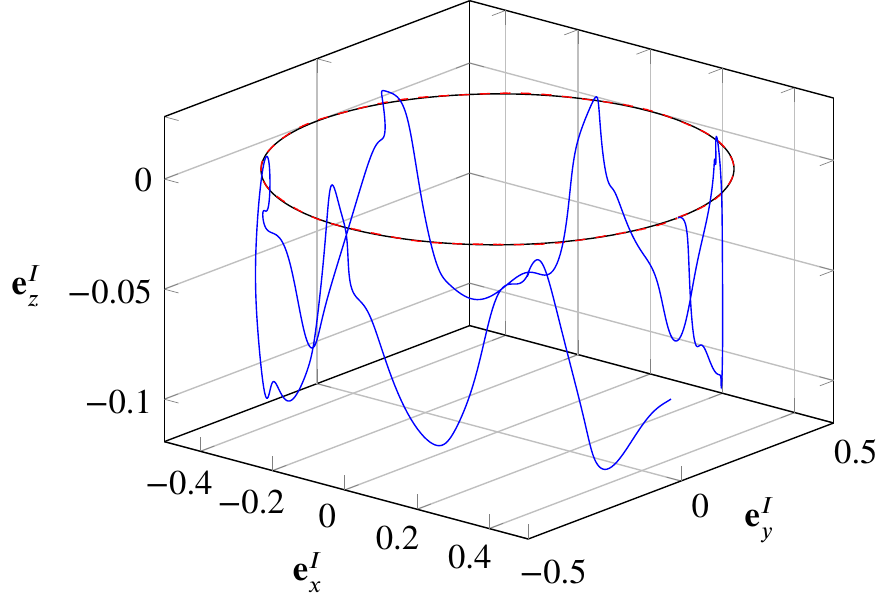}
	\end{subfigure}
	%
	\hspace{0.75cm}
	\begin{subfigure}[b]{0.47\textwidth}
		\caption{}
		\includegraphics{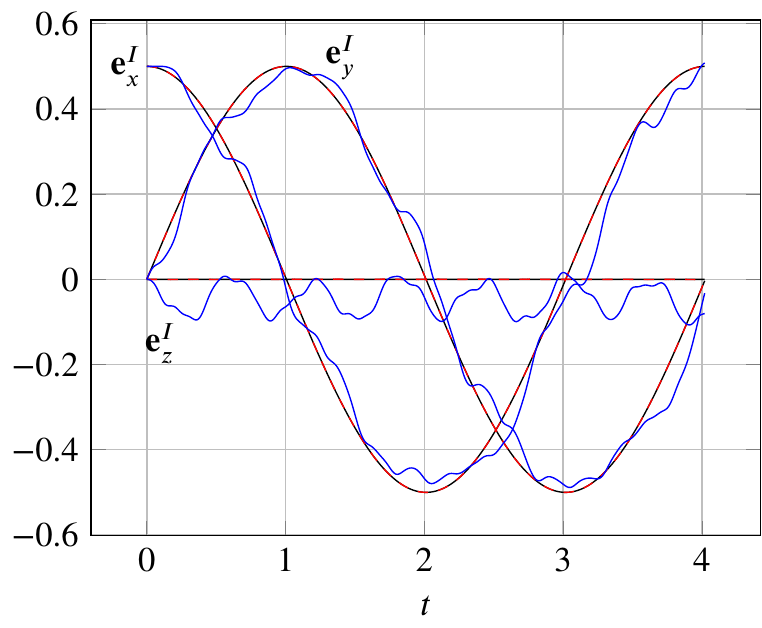}
	\end{subfigure}
	\caption{Tip displacement of the rigid top vs. stiff and soft rod solutions. The rigid top solution is drawn in black, the stiff rod in red and the soft rod in blue. (a) Spatial tip trajectory. (b) Tip displacements vs. time.}
	\label{fig:heavy_top}
\end{figure}

\begin{figure}[h!]
	\centering
	\hspace{0.2cm}
	\begin{subfigure}[b]{0.48\textwidth}
		\caption{}
		\includegraphics{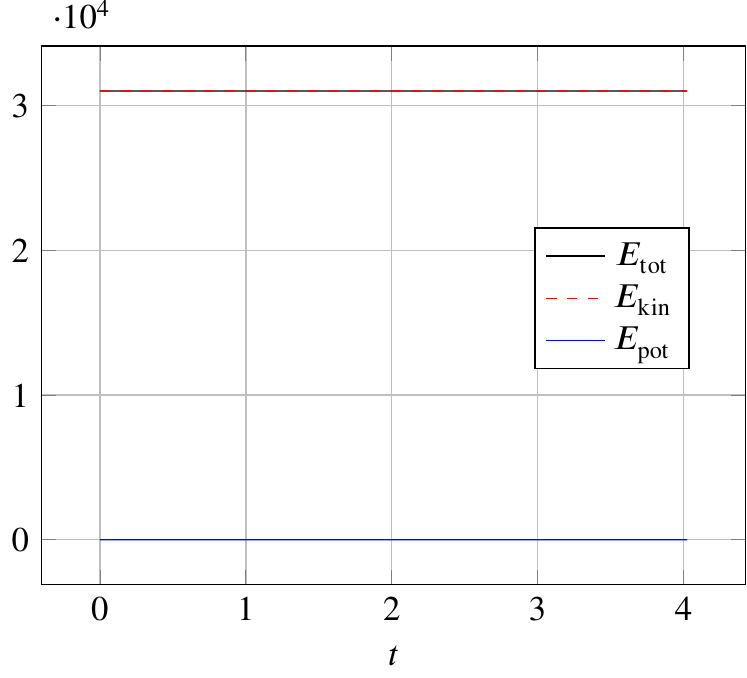}
	\end{subfigure}
	\hfill
	\begin{subfigure}[b]{0.48\textwidth}
		\caption{}
		\includegraphics{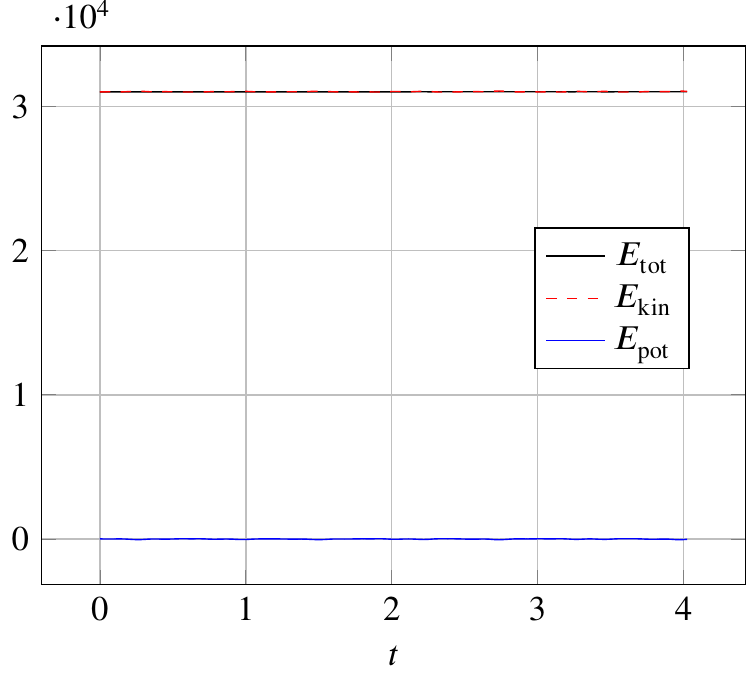}
	\end{subfigure}
	\caption{Time evolution of kinetic $E_\mathrm{kin}$, potential $E_\mathrm{pot}$ and total energy $E_\mathrm{tot}$ of the (a) stiff and (b) soft rod.}
	\label{fig:heavy_top_energy}
\end{figure}

In Figure~\ref{fig:heavy_top} (a), the spatial trajectories of the different tops' free ends are shown. When comparing the projections of the rod's free tip, it can be seen that the assertion is true that the solution of a stiff rod cannot be distinguished from the rigid body solution, while the soft rod's tip performs a fascinating oscillatory motion superimposed to the rigid body solution. The time evolution of potential $E_\mathrm{pot}$, kinetic $E_\mathrm{kin}$, and total energy $E_\mathrm{tot}$ for both the stiff and soft tops is illustrated in Figure~\ref{fig:heavy_top_energy}. Despite minor numerical artifacts, the total energy remains constant throughout the simulation, providing numerical evidence for the conservation of total energy as proven in of~\cite[Appendix D]{Harsch2023a}. However, since there is hardly any exchange of kinetic and potential energy in this example, the authors would like to emphasize that this is by no means a demonstrative study of the conservation properties for the presented family of rod finite element formulations.
\section{Conclusions}\label{sec:conclusions}

We presented an entire family of Petrov--Galerkin rod finite element formulations, where the individual formulations distinguish in the chosen interpolation strategy of centerline points and cross-section orientations, i.e., $\mR^{12}$-, $\mR^3 \times \SO(3)$-, or $\SE(3)$-interpolation. The Petrov--Galerkin projection method significantly simplifies the expressions of the discrete virtual work functionals, as no variations of the ansatz functions are required. Since all interpolation strategies are based on the nodal centerline points and cross-section orientations, they are independent of the underlying nodal $\SO(3)$-parameterization. We suggested a parametrization with nodal total rotation vectors whose singularity problems in the kinematic differential equations are resolved by introducing the complement rotation vector. Certainly, a singularity-free parametrization for instance with unit quaternions would be possible. However, an in-depth analysis of the effect of such a parametrization would have to be carried out. It is not clear whether the inclusion of the additional conditions that guarantee unit length of the quaternions influences the solvability of the equilibrium equations. Another interesting approach that follows from the formulation presented here would be to use constant strains as generalized position coordinates, cf. \cite{Renda2016}. The approximations of centerline and cross-section orientations would then follow from a recursive kinematic chain. Due to the Petrov--Galerkin method, the inertial virtual work would remain unaffected and the external virtual work would become only more complex if forces and moments are applied which are represented with respect to cross-section-fixed bases or the inertial basis, respectively.

While the $\SE(3)$-interpolation leads to an intrinsically locking-free formulation, the locking problems in the other interpolations could be cured by a simple reduced integration. The $\mR^{12}$-interpolation with polynomial degree $p=2$ stands out with a surprisingly low absolute error and a cubic spatial convergence behavior. Certainly, the comparison is not completely fair, as there would also be the possibility to interpolate the relative rotations \eqref{eq:R3xSO3_interpolation2} or twists \eqref{eq:SE3_interpolation2} by higher-order Lagrangian polynomials. This was already suggested by \cite{Crisfield1999} for the $\mR^3 \times \SO(3)$-interpolation. However in this case, the simplification to piecewise constant curvature \eqref{eq:discretized_picewise_curvature_measures} or strains \eqref{eq:discretized_picewise_strain_measures} is no longer valid and the evaluation of the $\SO(3)$ or $\SE(3)$-tangent maps, respectively, are required for computing the strains. This would make the already complex and computationally expensive interpolations even more demanding. In contrast, raising the polynomial order in the $\mR^{12}$-interpolation is not accompanied by an increased complexity in the formulation.

Consequently, for an application-oriented researcher, we recommend the $\mR^{12}$-interpolation strategy, which requires the least number of different concepts. These are for the static analysis (i) the continuous virtual work functionals of the Cosserat rod, (ii) the $\SO(3)$-parameterization with rotation vectors, (iii) $p$-th order Lagrangian polynomials, and for a standard Newton--Raphson algorithm (iv) the linearization of the equilibrium equations which can be found partially in \cite{Harsch2023a}.

\subsection*{Author contributions}
Both authors contributed equally to this publication.
%

\bibliographystyle{ieeetr}


\end{document}